\documentclass[notitlepage]{article} 
\usepackage{amsthm,amsmath,amssymb,latexsym,epsfig,wrapfig} 

\newlength{\jmr}
\newlength{\khov}
\newlength{\bernd}
\newlength{\milnor}
\newlength{\fw}
\newlength{\dima}
\newtheorem{rolle}{Rolle's Theorem} 

\newtheorem{kho}{Khovanski's Theorem on Real Fewnomials }

\newtheorem{nota}{Notation}

\newtheorem{lemma}{Lemma}
\newtheorem{prop}{Proposition}
\newtheorem{dfn}{Definition}

\newtheorem{main}{theorem} 
\newtheorem{thm}[main]{Theorem}
\newtheorem{des}{Univariate Generalized Descartes' Rule of Signs (UGDRS)}

\newtheorem{cor}{Corollary}
\newtheorem{rem}{Remark}	
\newtheorem{ex}{Example}

\newcommand{\vol}{\mathrm{Vol}} 
\newcommand{\inte}{\mathrm{Int}} 
\newcommand{\relint}{\mathrm{RelInt}}

\newcommand{\eps}{\varepsilon}

\newcommand{\cK}{\mathcal{K}}
\newcommand{\cN}{\mathcal{N}}

\newcommand{\supp}{\mathrm{Supp}}
\newcommand{\init}{\mathrm{Init}}
\newcommand{\conv}{\mathrm{Conv}}
\newcommand{\newt}{\mathrm{Newt}}

\newcommand{\thth}{{\underline{\mathrm{th}}}}
\newcommand{\rd}{ {\underline{ \mathrm{rd} } }  }
\newcommand{\st}{ {\underline{ \mathrm{st} } }  }

\newcommand{\Q}{\mathbb{Q}}
\newcommand{\R}{\mathbb{R}}
\newcommand{\C}{\mathbb{C}}
\newcommand{\N}{\mathbb{N}}
\newcommand{\Z}{\mathbb{Z}}

\newcommand{\Zn}{\Z^n}
 
\newcommand{\Rn}{\R^n}

\newcommand{\Rs}{\R^*}

\renewcommand{\qed}{$\blacksquare$}
\newcommand{\area}{\mathrm{Area}}
\newcommand{\dia}{$\diamond$}
\newcommand{\cL}{{\mathcal{L}}}

\newcommand{\cC}{\mathcal{C}}

\newcommand{\bO}{\mathbf{O}}

\begin{document}

\title{\mbox{}\\ 
\vspace{-1in}Counting Real Connected Components of Trinomial Curve 
Intersections and $m$-nomial Hypersurfaces}
% \subjclass{Primary 
% 34C08; % Connections with real algebraic geometry (fewnomials, 
%       % desingularization, zeros of Abelian integrals, etc.)
% Secondary 14P05, % Real algebraic sets 
% 30C15. % Zeros of polynomials, rational functions, and other 
                 % analytic functions (e.g. zeros of functions with bounded 
                 % Dirichlet integral) 
% } 

\author{\mbox{}\hspace{-.4in}\scalebox{.8}[1]{Tien-Yien Li\thanks{Li
was partially supported by a Guggenheim Fellowship and NSF Grant 
DMS-0104009.} 
\hspace{1.2in} J.\ Maurice Rojas\thanks{Rojas' work on this paper was partially 
supported by Hong Kong UGC Grant \#9040402-730, Hong Kong/France PROCORE 
Grant \#9050140-730, a grant from the Texas A\&M University Faculty of 
Science, and NSF Grant DMS-0211458.}  
\hspace{1.5in} Xiaoshen Wang\thanks{Some of Wang's research was done during a 
stay at the Key Laboratory for Symbolic Computation and Knowledge
Engineering of the Ministry of Education, P.\ R.\ China. Wang's research 
is supported in part by the Visiting Scholar
Foundation of Key Labs in Universities, Ministry of
Education, P.\ R.\ China. }} \\
\mbox{}\hspace{0in}\scalebox{.8}[1]{Department of Mathematics  
\hspace{.2in}
Department of Mathematics \hspace{.1in} Department of Mathematics and 
Statistics}\\
\mbox{}\hspace{0in}\scalebox{.8}[1]{Michigan State University\hspace{.6in} 
Texas A\&M 
University \hspace{.5in} University  of Arkansas at Little Rock}\\
\mbox{}\hspace{-.4in}\scalebox{.7}[1]{East Lansing, Michigan \ 48824 
\hspace{.2in} 
College Station, Texas \ 77843-3368  \hspace{.7in} Little Rock, Arkansas \ 
72204}\\
\mbox{}\hspace{-.5in}USA \hspace{1.5in} USA \hspace{1.8in} USA\\ 
\mbox{}\hspace{-.4in}{\tt li@math.msu.edu} \hspace{.3in} 
{\tt rojas@math.tamu.edu} \hspace{.7in} 
{\tt xxwang@ualr.edu }\\
\mbox{}\hspace{-.1in}\scalebox{.7}[1]{{\tt http://www.math.msu.edu/\~{}li }  
\hspace{.1in}{\tt http://www.math.tamu.edu/\~{}rojas}  \hspace{.2in} 
{\tt http://dms.ualr.edu/Faculty/Wang.html}} 
}

\date{\today} 

\maketitle 

\vspace{-.8cm}
\begin{center}
{\em In memory of Konstantin Alexandrovich Sevast'yanov, 1956--1984.} 
\end{center}

\begin{abstract} 
We prove that any pair of bivariate trinomials has at most 
$5$ isolated roots in the positive quadrant. The best previous upper 
bounds independent of the polynomial degrees were much larger, e.g., $248832$ 
(for just the non-degenerate roots) via a famous general result of Khovanski. 
Our bound is sharp, allows real exponents, allows degeneracies, and extends 
to certain systems of $n$-variate fewnomials, giving improvements over earlier 
bounds by a factor exponential in the number of monomials. We also derive 
analogous sharpened bounds on the number of connected components of the real 
zero set of a single $n$-variate $m$-nomial. 
\end{abstract} 

\section{Introduction} 
\label{sec:intro} 
Generalizing Descartes' Rule of Signs to multivariate systems of polynomial 
equations has 
proven to be a significant challenge. Recall that a weak version of this 
famous classical result asserts that any real univariate polynomial with 
exactly $m$ monomial terms has at most $m-1$ positive roots. This bound is 
sharp and generalizes easily to real exponents (cf.\ Section \ref{sec:back}). 
The original statement in Ren\'e Descartes' {\it La G\'eom\'etrie}\/ goes 
back to June of 1637, and Latham and Smith's English translation states that 
this result was observed even earlier by Thomas Harriot in his {\it
Artis Analyticae Praxis} (London, 1631) \cite[Footnote 196, 
Pg.\ 160]{descartes}. Proofs can be traced back to work of 
Gauss around 1828 and other authors earlier, but a definitive sharp bound for 
multivariate polynomial systems seems to have elluded us in the second 
millenium. This is particularly unfortunate since systems of sparse 
polynomial equations, and inequalities, now occur in applications as 
diverse as radar imaging \cite{keith}, chemistry \cite{gaterhub}, and neural 
net learning \cite{sagarrojas}. 

Here we take another step toward a sharp, higher-dimensional generalization 
of Descartes' bound by providing the first significant improvement on the 
case of curves in the plane, and certain higher-dimensional cases, since 
Khovanski's seminal work in the early 1980's \cite{kho}. Khovanski's 
revolutionary {\bf Theory of Fewnomials} \cite{kho,few} extends Descartes' 
bound to a broader class of analytic functions (incorporating certain measures 
of ``input complexity'') as well as higher dimensions, but the resulting 
bounds are impractically large even in the case of two variables. Our bounds 
are sharper than Khovanski's by a factor exponential in the 
number of monomial terms, allow degeneracies, and are optimal 
for the case of two bivariate trinomials. We then present similar sharpenings, 
also allowing real exponents and degeneracies, for the number of 
compact and non-compact connected components of the real zero set of a single 
sparse polynomial with any number of variables. 

\subsection{Main Results}\mbox{}\\ 
\label{sub:main} 
Perhaps the simplest way to generalize the setting of Descartes' Rule to 
higher dimensions and real exponents is the following: 
\begin{nota} 
For any $c\!\in\!\Rs\!:=\R\!\setminus\!\{0\}$ and 
$a\!=\!(a_1,\ldots,a_n)\!\in\!\Rn$, let 
$x^a\!:=\!x^{a_1}_1\cdots x^{a_n}_n$ and call $c x^a$ a 
{\bf monomial term}. We will refer to $\R^n_+\!:=\!\{x\!\in\!\Rn \; | 
\; x_i\!>\!0 \text{ \ for \ all \ } i \}$ as the {\bf positive orthant}. 
Henceforth, we will 
assume that $F\!:=\!(f_1,\ldots,f_k)$ where, for all $i$, 
$f_i\!\in\!\R[x^a \; | \; a\!\in\!\Rn]$ and $f_i$ has exactly 
$m_i$ monomial terms. We call $f_i$ an {\bf $\pmb{n}$-variate 
$\pmb{m_i}$-nomial} and, when $m_1,\ldots,m_k\!\geq\!1$, we call $F$ a 
{\bf $\pmb{k\times n}$ fewnomial system} (over $\R$) {\bf of type} 
$\pmb{(m_1,\ldots,m_k)}$. We call any homeomorphic image of the unit circle 
or a (closed, open, or half-open) interval an {\bf arc}.  
Finally, we say a real root $\zeta$ of $F$ is {\bf isolated} 
(resp.\ {\bf smooth}, {\bf non-degenerate}, or {\bf non-singular}) iff the only 
arc of real roots of $F$ containing $\zeta$ is $\zeta$ itself (resp.\ the 
Jacobian of $F$, evaluated at $\zeta$, has full rank). \dia 
\end{nota} 

\begin{dfn}
For any $m_1,\ldots,m_n\!\in\!\N$, let $\cN'(m_1,\ldots,m_n)$ 
(resp.\ $\cN(m_1,\ldots,m_n)$) denote the maximum 
number of non-degenerate (resp.\ isolated) roots an 
$n\times n$ fewnomial system of type $(m_1,\ldots,m_n)$ 
can have in the positive orthant. \dia
\end{dfn} 
Finding a tight upper bound on $\cN(m_1,\ldots,m_n)$ 
for $n\!\geq\!2$ remains a central problem in real algebraic geometry which 
is still poorly understood. For example, Anatoly Georievich Kushnirenko 
% Georgevich? 
conjectured 
in the mid-1970's \cite{kho} that $\cN'(m_1,\ldots,m_n)\!=\!\prod^n_{i=1} 
(m_i-1)$, at least for the case of integral exponents. The $n\times n$ 
polynomial system 
\begin{equation} 
\label{eqn:easy} 
\left(\prod^{m_1-1}_{i=1}(x_1-i),\ldots, \prod^{m_n-1}_{i=1}(x_n-i)\right) 
\end{equation}
was already known to provide an  
easy lower bound of $\cN'(m_1,\ldots,m_n)\!\geq\!\prod^n_{i=1} (m_i-1)$, 
but almost 30 years would pass until a counter-example to Kushnirenko's 
conjecture was published \cite{haas} (see Formula (\ref{eqn:haas}) in 
Section \ref{sub:rel}). However, it was known much earlier that Kushnirenko's 
conjectured upper bound could not be extended to $\cN(m_1,\ldots,m_n)$: The 
trivariate polynomial system
\begin{equation} 
\label{eqn:degen} 
\left(x(z-1),y(z-1),\prod^5_{i=1}(x-i)^2+\prod^5_{i=1}(y-i)^2 \right)
\end{equation} 
is of type $(2,2,21)$, has exactly $25$ ($>\!20\!=\!1\cdot 1\cdot 20$) 
roots in the positive 
octant, all of which are isolated and integral, but with Jacobian of rank 
$<\!3$ (see, e.g., \cite[note added in proof]{poly} or 
\cite[Ex.\ 13.6, Pg.\ 239]{ifulton}). Indeed, 
$\cN'(m_1,\ldots,m_n)\!\leq\!\cN(m_1,\ldots,m_n)$ for all 
$m_1,\ldots,m_n$, since a non-degenerate non-isolated root in $\Rn$ can never 
have more than $n-1$ tangent planes with linearly independent normal vectors. 
Cases where the inequality $\cN'(m_1,\ldots,m_n)\!\leq\!\cN(m_1,\ldots,m_n)$ 
is strict appear to be unknown.  

Interestingly, allowing degeneracies and real exponents introduces 
more flexibility than trouble in our approach: The proof of our first main 
result is surprisingly elementary, using little more than exponential 
coordinates and an extension of Rolle's Theorem from calculus. 
\begin{dfn}
\label{dfn:newt}
For any $S\!\subseteq\!\Rn$, let $\conv(S)$ denote the smallest
convex set containing $S$. Also, for any $m$-nomial of the form
$f(x)\!:=\!\sum_{a\in A} c_ax^a$, we call $\supp(f)\!:=\!\{a \; | \;
c_a\!\neq\!0\}$ the {\bf support} of $f$, and define
$\newt(f)\!:=\!\conv(\supp(f))$ to be the {\bf Newton polytope} of $f$. 
Finally, we let $Z_+(F)$ denote the zero set of $F$ in $\Rn_+$. \dia 
\end{dfn} 
\begin{thm} 
\label{thm:tri3} 
We have $\cN'(3,3)\!=\!\cN(3,3)\!=\!5$ and, more generally:  
\begin{itemize} 
\item[{\bf (a)}]{$\cN'(3,m)\!=\!\cN(3,m)\!\leq\!2^m-2$ for all 
$m\!\geq\!4$.} 
\item[{\bf (b)}]{ Any $n\times n$ fewnomial system $F\!:=\!(f_1,\ldots,f_n)$ 
of type $(m_1,\ldots,m_{n-1},m)$ with 
\begin{itemize}
\item[{\bf (i)}]{ $b_1+ \supp(f_1),\ldots,
b_{n-1}+\supp(f_{n-1})\!\subseteq\!A$ for some 
$b_1,\ldots,b_{n-1}\!\in\!\Rn$ and $A\!\subset\!\Rn$ of cardinality 
$n+1$} 
\item[{\bf (ii)}]{ $Z_+(f_1,\ldots,f_{n-1})$ smooth} 
\end{itemize} 
has no more than $n+n^2+\cdots+n^{m-1}$ isolated roots in $\Rn_+$, 
for all $m,n\!\geq\!1$. 
Also, for such $F$, the maximum number of 
non-degenerate and isolated roots in $\Rn_+$ are equal. }  
\item[{\bf (c)}]{ For any  
$\alpha_1,\alpha_2,\alpha_3,a_2,b_2,c_3,d_3,r_1,s_1,u_2,v_2\!\in\!\R$ 
and any degree $D$ polynomial $p\!\in\!\R[S_1,S_2]$ 
with $Z_+(p)$ smooth, 
the $2\times 2$ fewnomial system 
\[ (\star) \left\{ \ \begin{matrix}
\alpha_1+\alpha_2 x^{a_2}y^{b_2} + \alpha_3 x^{c_3}y^{d_3} \\ 
 p(x^{r_1}y^{s_1},x^{u_2}y^{v_2})  
\end{matrix}\right. 
\] 
has no more than $4\area(\newt(p))+2D+1$ ($\leq\!6D+1$) isolated roots 
in $\R^2_+$, where we normalize area so that the unit square 
has area $2$.  }  
\end{itemize} 
\end{thm} 
\noindent 
The quantities $\cN'(m_1,\ldots,m_n)$ and $\cN(m_1,\ldots,m_n)$ are much 
easier to compute when some $m_i$ is bounded above by $2$:  
all families of polynomial systems currently known to admit explicit formulae, 
including $n\times n$ {\bf binomial} systems, 
are summarized in Theorem \ref{thm:tri1} of Section \ref{sec:back}. 
\begin{rem} 
The value of $\cN(3,3)$ was previously unknown and there  
appears to be no earlier result directly implying the equality 
$\cN'(3,m)\!=\!\cN(3,m)$ for any $m$. In particular, the only other 
upper bound on $\cN'(3,m)$ or $\cN(3,m)$ until now was 
$\cN'(3,m)\!\leq\!3^{m+2}2^{(m+2)(m+1)/2}$, which evaluates to $248832$ 
when $m\!=\!3$. The best previous bound for the systems described 
in (b) and (c) were respectively $(n+1)^{m+n}2^{(m+n)(m+n-1)/2}$ and 
$1024D(D+2)^5$, counting only the non-degenerate roots. 
(See Khovanski's Theorem on Real Fewnomials in Section \ref{sub:rel} 
and Proposition \ref{prop:sandwich} of Section \ref{sec:back}.) \dia   
\end{rem} 
\begin{ex} 
Note that while we still don't know an upper bound on $\cN'(4,4)$ better 
than Khovanski's $4586471424$, we at least obtain a new approach for 
certain fewnomial systems with many monomial terms. For instance, 
part (c) of our first main theorem tells us that the $2\times 2$ 
fewnomial system: 
\[ \alpha_1+\alpha_2 x^{a_2}y^{b_2} + \alpha_3 x^{c_3}y^{d_3} \] 
\[ \beta_0 + \beta_{-1}x^{r_1}y^{s_1} 
+\beta_1x^{r_1+u_2}y^{s_1+v_2}+\beta_2x^{2(r_1+u_2)}y^{2(s_1+v_2)}+
\cdots+\beta_{100}x^{100(r_1+u_2)}y^{100(s_1+v_2)}\] 
has no more than $801\!=\!4\cdot 100+2\cdot 200+1$ isolated roots in $\R^2_+$, 
for all $\alpha_i,a_2,b_2,c_3,d_3,r_1,s_1,u_2,v_2\!\in\!\R$ and 
$\beta_i\!\in\!\R$ such 
that $\beta_0+\beta_1u+\cdots+\beta_{100}u^{100}$ has no degenerate roots. 
The pair of fewnomials $\left( (y-1)(y-2),y-\prod^{100}_{j=1}(x-j)\right)$ 
easily shows us that a system in this family can have as many as $200$ 
non-degenerate roots, if not more. Khovanski's Theorem on Real Fewnomials 
below yields an upper bound of $68878994643353600$ 
for just the number of non-degenerate roots. \dia 
\end{ex}  

We can also classify when a pair of bivariate trinomials 
has $5$ isolated roots in the positive quadrant via Newton polygons. 
In particular, note that while one can naturally associate a pair of polygons 
to $F$ when $n\!=\!2$, we can also associate a {\bf single} polygon by forming 
the {\bf Minkowski sum} $P_F\!:=\!\newt(f_1)+\newt(f_2)$. We can then give 
the following addendum to Theorem \ref{thm:tri3}.  
\begin{cor}
\label{cor:poly}  
A $2\times 2$ fewnomial system $F$ of type $(3,3)$ respectively has at most 
$0$, $2$, or $4$ isolated roots in $\R^2_+$, according as we restrict to those 
$F$ with $P_F$ a line segment, triangle, or $\ell$-gon with 
$\ell\!\in\!\{4,5\}$. 
\end{cor} 

The central observation that led to Theorem \ref{thm:tri3} may be of 
independent interest. 
We state it as assertion (5) of Theorem \ref{thm:cool} below. 
However, let us first define two more combinatorial quantities closely related 
to $\cN'$ and $\cN$. 
\begin{dfn}
For any $\mu,n\!\in\!\N$, we say that a $k\times n$ 
fewnomial system with exactly $\mu$ distinct exponent vectors 
is {\bf $\pmb{\mu}$-sparse}. Also, let $\cK'(n,\mu)$ (resp.\ 
$\cK(n,\mu)$) denote the maximum number of non-degenerate  
(resp.\ isolated) roots a $\mu$-sparse $n\times n$ fewnomial system 
can have in the positive orthant. \dia 
\end{dfn} 

Assertions (2) and (3) of our next main result dramatically refine the 
bounds of Oleinik, Petrovsky, Milnor, Thom, and Basu on the number of 
connected components of a real algebraic set \cite{op,milnor,thom,basu} in 
the special case of a single polynomial, and hold in the more general 
context of real exponents: 
\begin{thm} 
\label{thm:cool} 
Let $f$ be any $n$-variate $m$-nomial and let $P(n,m)$ denote be the 
maximum number of connected components of $Z_+(f)$ over all $n$-variate 
$m$-nomials. Also let 
$P_{\mathrm{comp}}(n,m)$ (resp.\ $P_{\mathrm{non}}(n,m)$)  
be the corresponding quantity counting just the compact (resp.\ non-compact) 
connected components. Finally, for any $r_1,s_1,u_2,v_2\!\in\!\R$ and 
any degree $D$ polynomial $p\!\in\!\R[S_1,S_2]$, let 
$\rho(x,y)\!:=\!p(x^{r_1}y^{s_1}, x^{u_2}y^{v_2})$. Then we have: 
\begin{enumerate} 
\addtocounter{enumi}{-1}
\item{ $P_{\mathrm{comp}}(n+1,2)\!=\!0$, $P_{\mathrm{non}}(n+1,2)\!=\!1$, 
$P_{\mathrm{comp}}(1,m)\!=\!m-1$, $P_{\mathrm{non}}(1,m)\!=\!0$,  
and $P_{\mathrm{non}}(n,0)\!=\!1$ for all $m,n\!\geq\!1$. } 
\item{$P_\mathrm{comp}(n,m)\!=\!0$ and 
$P_\mathrm{non}(n,m)\!=\!P_{\mathrm{non}}(m-1,m)\!\leq\!P(m-2,m)$ for 
$3\!\leq\!m\!\leq\!n+1$. } 
\item{$\max\left\{\left\lfloor\frac{m}{2}\right\rfloor-n-1,
\left(\left\lfloor\frac{m-1}{2n}\right\rfloor-1\right)^n\right\} 
\!\leq\!P_{\mathrm{comp}}(n,m)\!\leq\!2\lfloor\cK'(n,m)/2\rfloor$ for 
all $n\!\geq\!2$, and the last multiple of $2$ can be removed in the smooth 
case. Also, $Z_+(\rho)$ has no more than $\area(\newt(p))$ compact 
components, where we normalize area so that the unit square 
has area $2$. } 
\item{$\max\left\{m-1,\left(\left\lfloor\frac{m-1}{2(n-1)}\right\rfloor-1
\right)^{n-1}\right\}\!
\leq\!P_{\mathrm{non}}(n,m)\!\leq\!2P(n-1,m)$ for all $n\!\geq\!2$. 
Also, $Z_+(\rho)$ has no more than $2D$ non-compact connected components. } 
\end{enumerate} 
Furthermore, in the special case where $n\!=\!2$ and $Z_+(f)$ is 
smooth, let $I(m)$ (resp.\ $V(m)$) denote maximum number of 
isolated\footnote{ relative to the locus of points of inflection (resp.\
points of vertical tangency)...} 
inflection points (resp.\ isolated points of vertical tangency) 
of $Z_+(f)$. Then we also have 
\begin{enumerate}
\addtocounter{enumi}{3}
\item{$V(m)\!\leq\!\cK(2,m)$ for all $m\!\geq\!1$, and  
$Z_+(\rho)$ has no more than $\area(\newt(p))$ isolated points of vertical 
tangency. } 
\item{$I(m)\!\leq\!3\cK'(2,m)$ for all $m\!\leq\!3$, and 
$Z_+(\rho)$ has no more than $3\area(\newt(p))$ isolated inflection 
points.} 
\end{enumerate} 
In particular, $V(3)\!\leq\!1$ and $I(3)\!\leq\!3$, even if $Z_+(f)$ 
is not smooth. 
\end{thm}  
\noindent
Note that a non-compact connected component of $Z_+(f)$ can still have 
compact closure, since $\R^n_+$ is not closed in $\Rn$, e.g., 
$\left\{ (x_1,x_2)\!\in\!\R^2_+ \; | \; x^2_1+x^2_2\!=\!1\right\}$.  

While the above bounds on the number of connected components are non-explicit, 
they are stated so they can immediately incorporate any advance in 
computing $\cK'(n,\mu)$. So for a general and explicit upper 
bound independent of the underlying polynomial degrees now, one could, 
for instance, simply insert the explicit upper bound for $\cK'(n,\mu)$ 
appearing in Khovanski's Theorem on Real Fewnomials 
(see Section \ref{sub:rel} 
below) and the formula $\cK'(2,4)\!=\!5$ implied by 
Theorem \ref{thm:tri3} (see Proposition \ref{prop:sandwich} of Section 
\ref{sec:back}).  
\begin{cor}
\label{cor:exp} 
Following the notation of Theorem \ref{thm:cool}, 
$P(n,m)\!\leq\!\cK'(n,m)+2P(n-1,m)$ for all $n\!\geq\!2$. 
More explicitly, $P(n,m)\!\leq\!\sum^{n-1}_{i=0}2^i\cK(n-i,m)\!\leq\!
n(n+1)^m 2^{n-1}2^{m(m-1)/2}$. In particular, a tetranomial 
curve in $\R^2_+$ has no more than $4$ compact connected components and no 
more than\footnote{Theorem \ref{thm:cool} actually yields an upper bound 
of $6$ for the number of non-compact components so we cheated 
slightly by using Theorem \ref{thm:moment} below to get the very last bound.} 
$4$ non-compact connected components. \qed 
\end{cor} 
\noindent 
The bound above is already significantly sharper than an 
earlier bound of $(2n^2-n+1)^m (2n)^{n-1} 2^{m(m-1)/2}$, 
which held only for the smooth case, following from \cite[Sec.\ 3.14, 
Cor.\ 5]{few}. The bounds of Theorem \ref{thm:cool} 
also improve an earlier result of the middle author 
on smooth algebraic hypersurfaces \cite[Cor.\ 3.1]{real}. 

Our final main result shows us that we can considerably refine 
assertion (3) of Theorem \ref{thm:cool} if we take advantage of the 
underlying polyhedral structure. 
\begin{dfn} 
\label{dfn:newt2}
For any $w\!:=\!(w_1,\ldots,w_n)\!\in\!\Rn$ 
and any compact set $B\!\subset\!\Rn$, we let $B^w$ --- 
the {\bf face of $\pmb{B}$ with
inner normal $\pmb{w}$} --- be the set of all $x\!\in\!B$
minimizing the inner product $w\cdot x$.
Finally, for any any $n$-variate $m$-nomial $f$ of the form 
$\sum_{a\in A} c_ax^a$, we let $\init_w(f)$
--- the {\bf initial terms} of $f$ (with respect to $w$) ---
be $\sum_{a\in A^w} c_a x^a$. \dia 
\end{dfn} 

\noindent
Recall that the dimension of a polytope $P\!\subseteq\!\Rn$ is the dimension 
of the smallest subspace containing a translate of $P$ and that a {\bf facet} 
of an $n$-dimensional polytope is simply a face of dimension $n-1$. 
\begin{thm}
\label{thm:moment}
Let $f$ be any $n$-variate $m$-nomial $f$ with
$n$-dimensional Newton polytope. Assume further that $Z_+(\init_w(f))$
is smooth for all $w\!\in\!\Rn\!\setminus\!\{\bO\}$.
Then the number of non-compact connected components of $Z_+(f)$ is no more
than \[ \sum\limits_{\substack{w \text{ a unit inner facet } \\ 
\text{ normal of } \newt(f)}} \!\!\!\!\!\!\!\!\!\!\!\!N_w, \] where $N_w$ 
denotes the number of the number of connected components of $Z_+(\init_w(f))$. 
In particular,  this bound is no larger than $\!\!\!\!\!\!\sum\limits_{Q 
\text{ a facet of } P} \!\!\!\!\!\!P(n-1,\#\supp(f)\cap Q)$. Finally, when 
$n\!=\!2$ and $Z_+(f)$ is smooth as well, the last upper bound is sharp 
and can be simplified to $\lfloor m'/2\rfloor$, where $m'$ is the number
of points of $\supp(f)$ lying on the boundary of $\newt(f)$.
\end{thm}

\noindent
Note that $Z_+(f)$ need not be smooth and our bound above is 
completely independent of the number of exponent vectors lying in 
the interior of $\newt(f)$. The bivariate example $f(x,y)\!=\!
y-\prod^{m'-2}_{i=1} (x-i)$ easily shows that the very last 
bound is sharp. A more intricate trivariate example follows. 
\begin{ex}
Taking $n\!=\!3$, suppose $f$ is
\[ \alpha_1+\alpha_2x^{3a}+\alpha_3z^{3c}+\alpha_4x^{3a}z^{3c} +
\beta_1x^a y^b z^c+\beta_2 x^{2a}y^b z^c + \beta_3 x^a y^b z^{2c}
+ \beta_4 x^{2c} y^b z^{2c} + \sum^K_{i=1} \gamma_i x^{a_i}y^{b_i}
z^{c_i},\]
where $K$ is any positive integer, the $\alpha_j$, $\beta_j$, $\gamma_j$
are any nonzero real constants, $a,b,c\!>\!0$ and, for all $i$,
$a\!<\!a_i\!<\!2a$, $c\!<\!c_i\!<\!2c$, and $0\!<\!b_i\!<\!b$. Note that
$\newt(f)$ is a snub pyramid with a rectangular base and thus has the same
face lattice as a cube. Note also that no exponent vector of $f$ lies in the
relative interior of any face of $\newt(f)$ of dimension $1$ or $2$. It is 
then easily checked that
\[ (\alpha_1\alpha^2_4-\alpha_2\alpha_3\alpha_4)
(\beta_1\beta^2_4-\beta_2\beta_3\beta_4)
Q(\alpha,\beta)\neq 0,\]

\noindent
where $Q$ is a product of $4$ more complicated polynomials, is a
sufficient condition for all the $Z_+(\init_w(f))$ to be smooth.
So, under the last assumption,
Theorem \ref{thm:moment} tells us that the zero set of $f$ in the
positive octant has no more than $6P(2,4)\!\leq\!6\cdot
(4+6)\!=\!60$ non-compact connected components, employing corollary
\ref{cor:exp} and the obvious fact that
$P(n,m)\!\leq\!P_\mathrm{comp}(n,m)+P_\mathrm{non}(n,m)$ for the first
inequality. Note that Theorem \ref{thm:cool} would have given us
a less explicit upper bound of $P_\mathrm{non}(3,8+K)$ which, by
assertion (3), exceeds $60$ for all $K\!\geq\!37$ (if not earlier). \dia
\end{ex}

Note also that the assumption on the $\init_w(f)$ is rather mild:
it follows easily from Sard's Theorem \cite{hirsch} that
our smoothness condition will hold for
a generic choice of the coefficients of $f$, e.g., all coefficient vectors
outside a set of measure zero in $\C^{\#\supp(f)}$ depending only on
$\supp(f)$. In particular, this hypothesis can become vacuous depending 
on the underlying Newton polytope. 
\begin{cor} 
\label{cor:overdet} 
Following the notation above, assume that $\newt(f)$ is 
{\bf simplicial} (i.e., for all $d\!<\!\dim \newt(f)$ every $d$-dimensional 
face of $\newt(f)$ has exactly $d+1$ vertices) and that [the relative 
interior of a face $Q$ of $\newt(f)$ contains a point of 
$\supp(f) \Longrightarrow Q$ is a vertex]. Then $Z_+(\init_w(f))$ is smooth 
for all $w\!\in\!\Rn\!\setminus\!\{\bO\}$. \qed 
\end{cor} 
\noindent 
Corollary \ref{cor:overdet} follows easily 
from the fact that for such an $f$, and any vector 
$w\!\in\!\Rn\!\setminus\!\{\bO\}$, 
$Z_+(\init_w(f))$ is analytically diffeomorphic to $\R^{n-1}_+$. The 
latter fact in turn follows easily via a monomial change of variables 
(cf.\ Proposition \ref{prop:easy} of Section \ref{sec:back}). 

\subsection{Important Related Results}\mbox{}\\ 
\label{sub:rel}  
The only available results for bounding the number of real roots, 
other than those coming from Fewnomial Theory \cite{bocook,grigoriev,risler,
few,zell,real}, depend strongly on the individual exponents of $F$ 
and are actually geared more toward counting complex roots, e.g., 
\cite{bkk,kaza,blr,jpaa}. (We also note that while the 
bounds of \cite{zell} generalize Khovanski's theory to solution 
sets of {\bf inequalities} involving Pfaffian functions, they do 
not appear to yield any new bounds on the quantities $\cK$ and $\cN$ 
we study.) So proving just $\cN(3,3)\!<\!\infty$ already suggests an  
analytic approach. Nevertheless, the bounds from \cite{bkk,kaza,blr,jpaa} 
can be quite practical when the exponents are integral and the degrees 
of the polynomials are small. 

Let us also point out that the term ``fewnomial'' is due to 
Kushnirenko and that the first explicit bounds 
in Fewnomial Theory were derived (not yet in complete generality) 
by Konstantin Alexandrovich Sevast'yanov in unpublished work around 
1979 \cite{askold}. Dima Yu.\ Grigoriev and Askold Georgevich Khovanski 
have also pointed out that shortly after Kushnirenko formulated his 
conjecture, a simple counter-example with $n\!=\!2$ was found by 
a student at Moscow State University \cite{dima,askold}. Unfortunately, while 
the counter-example was verified by Khovanski himself \cite{askold}, it does 
not seem to have been recorded and the name of its inventor (who 
left mathematics immediately after graduating) seems to have been forgotten. 

As for the size of upper bounds on the number of real roots, 
it is interesting to note that the best current general bounds independent 
of the polynomial degrees are exponential in the number of monomial terms of 
$F$, even for fixed $n$. Observe one of the masterpieces of real algebraic 
geometry. 
\begin{kho}
(See \cite{kho} and \cite[Cor.\ 6, Pg.\ 80, Sec.\ 3.12]{few}.)  
We have $\cK'(n,\mu)\!\leq\!(n+1)^\mu 2^{\mu(\mu-1)/2}$. 
More generally, the $n\times n$ fewnomial system
\[ q_1(x)=\cdots =q_n(x)=0,\] where each $q_j$ is a 
polynomial of degree 
$D_i$ in $x_1,\ldots,x_n$ and $x^{a_1},\ldots,x^{a_\mu}$ for 
some $a_1,\ldots,a_\mu\!\in\!\Rn$, has no more than 
$2^{\mu(\mu-1)/2}(1+\sum^n_{i=1} D_i)^\mu\prod^n_{i=1} D_i$ non-degenerate 
roots in $\Rn_+$. \qed 
\end{kho} 

Finding non-trivial lower bounds on even $\cK'(2,\mu)$ seems quite hard and 
surprisingly little else is known about what an optimal version of Khovanski's 
Theorem on Real Fewnomials should resemble. For example, around 
1996, Ilya Itenberg and Marie-Fran\c{c}oise Roy proposed 
a conjectural polyhedral generalization of Descartes' 
Rule to multivariate systems of equations \cite{ir}, based on a famous 
construction from Oleg Viro's 1983 Leningrad thesis (see, e.g., 
\cite{viro}) and later extensions by Bernd Sturmfels 
\cite{berndviro}. A bit later, Sturmfels offered US\$500 
for a proof that Itenberg and Roy's proposed upper bound held for 
the following family of $2\times 2$ systems of type $(4,4)$:  
\begin{equation}  
\label{eqn:bernd} 
(-x^5+a_1y^5+a_2x^3y^5+a_3x^6y^8, 
-y^5+b_1x^5+b_2x^5y^3+b_3x^8y^6),  
\end{equation} 
where $a_1,a_2,a_3,b_1,b_2,b_3\!>\!0$. The Itenberg-Roy conjecture yields 
an alleged upper bound of $3$ for this family, and Jeff Lagarias and Thomas 
Richardson later won Sturmfels' prize by showing that this bound in 
fact holds \cite{lr}. However, the Itenberg-Roy conjecture was later 
invalidated by the $2\times 2$ system 
\begin{equation}
\label{eqn:liwang} 
(y-x-1,y^3+0.01x^3y^3-9x^3-2) 
\end{equation} 
found by the left and right authors (Li and Wang): 
this system has exactly $3$ roots in $\R^2_+$, whereas the conjectured bound 
would have only been $2$ \cite{liwang}. 
Perhaps the most important counter-example in this growing theory is 
Haas' recent counter-example to Kushnirenko's Conjecture: It is 
\begin{equation} 
\label{eqn:haas} 
(x^{108}_1+1.1x^{54}_2-1.1x_2, x^{108}_2+1.1x^{54}_1-1.1x_1),
\end{equation} 
which has $5$ ($>\!4\!=\!2\cdot 2$) roots in the positive quadrant
\cite{haas}. Jan Verschelde has also verified numerically via his
software package {\tt PHCPACK} \cite{verschelde} that
there are exactly $108^2\!=\!11664$ complex roots, and thus (assuming the
floating-point  calculations were sufficiently good) each root is
non-degenerate by B\'ezout's theorem on intersections of complex hypersurfaces
\cite[ex.\ 1, pg.\ 198]{sha}. 

As for asymptotic behavior, it is still unknown whether even $\cK'(2,\mu)$ is 
polynomial in $\mu$: even the special case of $2\times 2$ fewnomial 
systems of type $(3,m)$ is still open. Note also that this kind of 
polynomiality requires the number of variables to be fixed: 
the system \mbox{$(x^2_1-3x_1+2,$} $\ldots,x^2_n-3x_n+2)$ shows us that 
$\cN'(\underset{n}{\underbrace{3,\ldots,3}})$ is already \mbox{exponential in 
$n$.} More to the point, it is also unknown whether a simple modification 
of Kushnirenko's conjectured bound 
(e.g., increasing the formula by a constant power or a factor exponential in 
$n$) would at last yield a true, sharp, and general improvement of 
Khovanski's Theorem on Real Fewnomials. The $2k\times 2k$ fewnomial system 
\[ (x^{108}_1+1.1y^{54}_1-1.1y_1,y^{108}_1+1.1x^{54}_1-1.1x_1,\ldots,
x^{108}_k+1.1y^{54}_k-1.1y_k,y^{108}_k+1.1x^{54}_k-1.1x_k),\] thanks to 
Haas' counter-example, easily shows that one  
needs at least an extra multiple no smaller than 
$\left(\frac{\sqrt{5}}{2}\right)^n$ if some Kushnirenko-like bound is to be 
salvaged.  

Another question with even deeper implications is whether 
there is an algorithm for approximating the real roots of a fewnomial system 
whose complexity depends mainly on the number of {\bf real} roots. Since all 
current algorithms for real-solving have complexity 
bounds essentially matching the analogous bounds for solving over 
the complex numbers, a positive answer would yield tremendous 
speed-ups, both practical and theoretical, for real-solving. However, little 
is known beyond the special cases of $n\times n$ binomial systems 
\cite[Main Thm.\ 1.3]{real} and univariate polynomials with $3$ monomial terms 
or less \cite{rojasye}: For these cases, one can indeed obtain algorithms 
beating the known lower bounds \cite{renegar} for solving over the complex 
numbers, and \cite{rojasye} also shows that one can at least find the 
isolated inflection points and vertical tangents of a trinomial 
curve about as quickly. 

Let us conclude our introduction with a recent number-theoretic parallel: It 
has just been shown by the middle author \cite{ari} that the number of 
{\bf geometrically} isolated\footnote{A root is geometrically 
isolated iff it is a zero-dimensional component of the underlying 
zero set in $\bar{\cL}^n$, where $\bar{\cL}$ is the algebraic closure of 
$\cL$.} roots in $\cL^n$ of any $\mu$-sparse $k\times n$ 
polynomial system, {\bf over any $\pmb{\mathfrak{p}}$-adic field $\pmb{\cL}$}, 
is no more than $1+\left(\cC_\cL n(\mu-n)^3\log \mu \right)^n$, where 
$\cC_\cL$ is a constant depending only on $\cL$ (see also \cite{myadic} and 
the references therein for earlier results in this direction). In particular, 
since $\Q\!\subset\!\Q_2$, one thus obtains a bound on the number of isolated 
roots in $\cL^n$ which is polynomial in $\mu$ for fixed $n$, with $\cL$ now 
any fixed number field. One should note that 
$\mathfrak{p}$-adic 
fields, just like $\R$, are complete with respect to a suitable metric. So 
there appears to be a deeper property of metrically complete fields lurking 
in these quantitative results. 
\begin{rem} 
\label{rem:napo} 
Domenico Napoletani has recently shown that to calculate 
$\cN'(m_1,\ldots,m_n)$ for any given $(m_1,\ldots,m_n)$, it suffices 
to restrict to the case of integral exponents \cite{napo}. Here, we 
will bound $\cN(n+1,\ldots,n+1,m)$ directly, in the aforementioned 
cases, without using this reduction. \dia 
\end{rem} 

\subsection{Organization of the Proofs and Obstructions to Extensions}\mbox{}\\ 
Section \ref{sec:back} provides some background and unites some simple cases 
where Kushnirenko's conjectured bound in fact holds, and the equalities 
$\cK'(n,\mu)\!=\!\cK(n,\mu)$ and 
$\cN'(m_1,\ldots,m_n)\!=\!\cN(m_1,\ldots,m_n)$ are true. We then prove 
Theorem \ref{thm:tri3} in Sections \ref{sec:tx} and \ref{sec:tri}, and 
prove Theorem \ref{thm:cool} in Section \ref{sec:morse}. 
Proving the (restricted) upper bound on $\cN(n+1,\ldots,n+1,m)$ turns out to 
be surprisingly elementary, but lowering the upper bound on $\cN(3,3)$ to $5$ 
then becomes a more involved case by case analysis. Section \ref{sec:moment} 
then applies a variant of the momentum map from symplectic/toric geometry 
(see, e.g., \cite{smale,souriau} and \cite[Sec.\ 4.2]{tfulton}) 
to prove Theorem \ref{thm:moment}.  

Section \ref{sec:tri} gives an alternative geometric proof that 
$\cN(3,3)\!\leq\!6$. We include this second proof for motivational purposes 
since it appears to be the first known improvement over 
$\cN'(3,3)\!\leq\!248832$, and since it is the only approach 
we know which yields part (c) of Theorem \ref{thm:tri3}. 

The reader should at this point be aware that our results 
can of course be combined and interweaved to generate 
much more complicated examples (with more monomial terms,  
more complicated supports, and more variables) which admit upper 
bounds on the number of roots in $\R^n_+$ significantly sharper than 
Khovanski's Theorem on Real Fewnomials (see, e.g., Theorem \ref{thm:tri1} 
and the paragraph after in the next section). Nevertheless, it should also 
be clear that there are still many simple fewnomial systems where nothing 
better than Khovanski's bound is available, e.g., the exact values of 
$\cN'(4,4)$ and $\cN(4,4)$ remain unknown. So let us close with some brief 
remarks on the obstructions to extending Theorem \ref{thm:tri3} to more 
complicated fewnomial systems. In particular, the two main techniques we use 
are (A) a recursion involving derivatives of certain analytic functions, and 
(B) an extension of Rolle's Theorem (cf.\ Section \ref{sec:back}) to 
intersections of lines with certain fewnomial curves. 

Our technique from (A) succeeds precisely because the 
underlying recursion stops in a number of steps depending only on $m$ 
and $n$. In particular, while one can apply the same technique to certain 
slightly more complicated systems (cf.\ the proof of part (b) of Theorem 
\ref{thm:tri3} in Section \ref{sec:tx}), applying the same technique to a 
system of type $(4,m)$ results in a much more complicated recursion which 
won't terminate without strong restrictions on the exponents; 
and even then the number of steps begins to depend on the exponents.  
The geometric reason for this is that we in essence project 
our roots to a line to start our recursion, and such 
projected roots appear to satisfy sufficiently simple equations 
just for the systems defined in part (b) (see Remark \ref{rem:proj} of 
the next section). 

Our technique from (B) succeeds for the systems $(f_1,f_2)$ coming from 
($\star$) precisely because (i) $Z_+(f_1)$ is diffeomorphic to 
a line in a very special way, and (ii) the equations arising from checking 
inflection points and vertical tangents of $Z_+(f_2)$ have a fewnomial 
structure very similar to that of $f_2$. In particular, for the systems in 
($\star$), we construct our stated bound by a special application of 
Bernstein's Theorem \cite{bkk} in the $2\times 2$ case. However, increasing 
the number of variables of $p$ (i.e., the number of monomial terms of $f_2$) 
leaves us with a system of equations apparently not reducible to Bernstein's 
Theorem.  

Nevertheless, we suspect that there are many similar 
improvements to Fewnomial Theory over $\R$ which are quite tractable, and 
we hope that our paper serves to inspire more activity in this area. 

\section{The Pyramidal, Simplicial, and Zero Mixed Volume Cases} 
\label{sec:back} 
Let us first note some simple inequalities relating the quantities 
$\cK'$, $\cK$, $\cN'$, and $\cN$. 
\begin{prop}
\label{prop:sandwich}
We have $(\mu-1)^n\!\leq\!\cK'(n,\mu)\!\leq\!\cK(n,\mu)$,\\
\scalebox{.9}[1]{$\cN'(m_1,\ldots,m_n)\!\leq\!\cK'(n,
m_1+\cdots+m_n-n+1)\!\leq\!\cN'(\underset{n}
{\underbrace{m_1+\cdots+m_n-2n+2,\ldots,m_1+\cdots+m_n-2n+2}})$,}\\ 
and\\ 
\scalebox{.9}[1]{$\cN(m_1,\ldots,m_n)\!\leq\!\cK(n,
m_1+\cdots+m_n-n+1)\!\leq\!\cN(
\underset{n}{\underbrace{m_1+\cdots+m_n-2n+2,\ldots,m_1+\cdots+m_n-2n+2}})$
,}\\
where we set $\cN(m_1,\ldots,m_n)\!=\!\cN'(m_1,\ldots,m_n)\!=\!0$ if any
$m_i$ is negative. In particular, by Theorem \ref{thm:tri3}, we thus
have $\cK'(2,4)\!=\!\cK(2,4)\!=\!5$. \qed
\end{prop}
 
\noindent
Indeed, the last two ``left-hand'' inequalities follow
simply by dividing each $f_i$ by a suitable monomial, while
Gaussian elimination on the monomial terms of $F$ easily yields
the last two ``right-hand'' inequalities.

Let us next give a simple geometric characterization of certain 
fewnomial systems that admit easy root counts. 
\begin{dfn}
Let us call any collection $L_1\!\subsetneqq\cdots
\subsetneqq\!L_n\!=\!\Rn$ of $n$ non-empty subspaces of $\Rn$ 
(so that $\dim L_i\!=\!i$ for all $i$) a {\bf complete flag}. 
Noting that any polytope in $\Rn$
naturally generates a subspace 
of $\Rn$ via the set of linear
combinations of all {\bf differences} of its vertices, let 
$F\!=\!(f_1,\ldots,f_n)$ be an $n\times n$ fewnomial system and, for all
$i$, let $L_i$ be the linear subspace so generated by $\newt(f_i)$. 
We then say that $F$ is {\bf pyramidal} iff the Newton polytopes of 
$F$ generate a complete flag. Finally, letting 
$A\!:=\![a_{ij}]$ be any real $n\times n$ matrix, $x\!:=\!(x_1,\ldots,x_n)$, 
$y\!:=\!(y_1,\ldots,y_n)$, and $y^A\!:=\!(y^{a_{11}}_1\cdots
y^{a_{n1}}_n,\ldots, y^{a_{1n}}_1\cdots y^{a_{nn}}_n)$, we call any 
change of variables of the form 
$x\!=\!y^A$ a {\bf monomial change of variables}.  \dia 
\end{dfn}

\noindent
For example, the systems from (\ref{eqn:easy}) (cf.\ Section \ref{sub:main}) 
are pyramidal, but systems (\ref{eqn:degen}) (cf.\ Section \ref{sub:main}), 
(\ref{eqn:bernd}), (\ref{eqn:liwang}), and (\ref{eqn:haas})  
(cf.\ Section \ref{sub:rel}) are all non-pyramidal. Note in particular that 
all binomial systems are pyramidal, but a $2\times 2$ fewnomial system of 
type $(3,3)$ certainly need not be pyramidal. Pyramidal systems are a 
simple generalization 
of the so-called ``triangular'' systems popular in Gr\"obner-basis papers on
computer algebra. The latter family of systems simply consists of those $F$
for which the variables can be reordered so that 
for all $i$, $f_i$ depends only on $x_1,\ldots,x_i$. Put another 
way, pyramidal systems are simply the image of a triangular system (with {\bf 
real} exponents allowed) after multiplying the individual equations by 
arbitrary monomials and then performing a monomial 
change of variables. 

Recall that an {\bf analytic} subset of a domain $U\!\subseteq\!\Rn$ is simply 
the zero set of an analytic function defined on $U$. We then have the 
following fact on monomial changes of variables. 
\begin{prop}
\label{prop:easy} 
If $x\!:=\!(x_1,\ldots,x_n)\!\in\!\Rn_+$ and $A$ is a real invertible  
$n\times n$ matrix, then 
$(x^{A})^{A^{-1}}\!=\!x$ and the {\bf monomial map} defined by $x\mapsto x^A$ 
is an analytic automorphism of the positive orthant. In particular, 
such a map preserves smooth points, singular points, and the 
number of compact and non-compact connected components, of 
analytic subsets of the positive orthant. 
Furthermore,  this invariance also holds for real $m$-nomial zero sets in the 
positive orthant. \qed 
\end{prop} 
\noindent 
The assertion on analytic subsets follows easily from an application of 
the chain rule from calculus, and noting that such monomial maps 
are also diffeomorphisms. That the same invariance holds for $m$-nomial 
zero sets follows immediately upon observing that the substitution 
$(x_1,\ldots,x_n)\!=\!(e^{z_1},\ldots,e^{z_n})$ maps 
any $n$-variate real $m$-nomial to a real analytic function, and noting that 
the map defined by $(t_1,\ldots,t_n)\mapsto 
(e^{t_1},\ldots,e^{t_n})$ is a diffeomorphism from $\Rn$ to $\Rn_+$. 
\begin{rem}
The real zero set of $x_1+x_2-1$, and the change of variables 
$(x_1,x_2)\!=\!\left(\frac{y_1}{y_2},y_1y_2\right)$, show that 
the number of isolated {\bf inflection} 
points need not be preserved by such a map: the underlying curve 
goes from having no isolated inflection points  
to having one in the positive quadrant. \dia 
\end{rem} 

We will later need the following analogous geometric extension of 
the concept of an over-determined system.  
\begin{dfn} 
\label{dfn:zero} 
Given polytopes $P_1,\ldots,P_n\!\subset\!\Rn$, we say that they have {\bf 
mixed volume zero} iff for some $d\!\in\!\{0,\ldots,n-1\}$ there exists a 
$d$-dimensional subspace of $\Rn$ containing translates of $P_i$ for at least 
$d+1$ distinct $i$. \dia    
\end{dfn} 
\noindent 
The {\bf mixed volume}, originally defined by Hermann Minkowski in the 
late 19$^\thth$ century, is a nonnegative function defined for 
all $n$-tuples of convex bodies in $\Rn$, and satisfies many 
natural properties extending the usual $n$-volume. 
The reader curious about mixed volumes of polytopes in 
the context of solving polynomial equations can consult \cite{buza,jpaa} 
(and the references therein) for further discussion. A simple special case 
of an $n$-tuple of polytopes with mixed volume zero is the $n$-tuple of 
Newton polytopes of an $n\times n$ fewnomial system 
where, say, the variable $x_i$ does not appear. By multiplying 
the individual $m$-nomials by suitable monomials, and applying a 
suitable monomial change of variables, the following 
corollary of Proposition \ref{prop:easy} is immediate. 
\begin{cor} 
\label{cor:zero} 
Suppose $F$ is a fewnomial system, with only finitely many roots 
in the positive orthant, whose $n$-tuple of Newton polytopes has mixed volume 
zero. Then $F$ has no roots in the positive orthant. \qed 
\end{cor} 
\noindent 
Indeed, modulo a suitable monomial change of variables, one need only observe 
that the existence of a single root in the positive orthant implies the 
existence of an entire ray of roots (parallel to some coordinate axis) in the 
positive orthant. 

We will also need the following elegant extension of Descartes' Rule to real 
exponents. It's proof involves a very simple induction using Rolle's 
Theorem (cf.\ the next section) and dividing by suitable monomials 
\cite{few} --- tricks we will build upon in the next section. 
\begin{dfn}
For any sequence $(c_1,\ldots,c_m)\!\in\!\R^m$, its {\bf number of 
sign alternations} is the number of pairs $\{j,j'\}\!\in\{1,\ldots,m\}$ such 
that $j\!<\!j'$, $c_jc_{j'}\!<\!0$, and $c_i\!=\!0$ when $j\!<\!i\!<\!j'$. \dia 
\end{dfn} 
\begin{des}
\label{thm:des} 
Let $c_1,a_1,\ldots,c_m,a_m$ be any real 
numbers with $a_1\!<\cdots <\!a_m$. 
Then the number of positive roots of $\sum^m_{i=1}c_ix^{a_i}_1$ is 
at most the number of sign alternations in the sequence 
$(c_1,\ldots,c_m)$. In particular, 
$\cK'(1,m)\!=\!\cK(1,m)\!=\!\cN'(m)\!=\!\cN(m)\!=\!m-1$. \qed 
\end{des} 

As a warm-up, we can now prove a stronger version of Kushnirenko's 
conjecture for certain fundamental families of special cases. In particular, 
we point out that aside from the domain of Theorem \ref{thm:tri3}, the 
equalities $\cK'(n,\mu)\!=\!\cK(n,\mu)$ and 
$\cN'(m_1,\ldots,m_n)\!=\!\cN(m_1,\ldots,m_n)$ appear to be known only for 
the cases stated in UGDRS and assertions (0), (2), and (4) below. 
\begin{thm}
\label{thm:tri1} 
Suppose $F$ is an $n\times n$ fewnomial system of type 
$(m_1,\ldots,m_n)$ (so $m_1,\ldots,m_n\!\geq\!1$) and consider  
the following independent conditions: 
\begin{itemize} 
\item[{\bf (a)}]{The $n$-tuple of Newton polytopes of $F$ 
has mixed volume zero.} 
\item[{\bf (b)}]{All the supports of $F$ can be translated into 
a single set of cardinality $\leq\!n+1$.}  
\item[{\bf (c)}]{$F$ is pyramidal.} 
\end{itemize} 
Then, following the notation of Theorem \ref{thm:tri3}, we have: 
\begin{enumerate} 
\addtocounter{enumi}{-1} 
\item{$\cN(m_1,\ldots,m_n)$ is $0$, $1$, or 
$\prod^n_{i=1}(m_i-1)$ if we respectively restrict to case (a), (b), or (c). 
Also, in all these cases, $\cN'(m_1,\ldots,m_n)\!=\!\cN(m_1,\ldots,m_n)$.} 
\item{If (a), (b), or (c) hold then [$F$ has infinitely many roots 
$\Longrightarrow F$ has {\bf no} isolated roots].}
\item{$\cN'(m_1,m_2,\ldots,m_n)\!=\!\cN(m_1,m_2,\ldots,m_n)\!=\!0 
\Longleftrightarrow$ some $m_i$ is $\leq\!1$ } 
\item{$m_1\!=\!2 \Longrightarrow$  
[$\cN(m_1,m_2,\ldots,m_n)\!=\!\cN(m_2,\ldots,m_n)$ and 
$\cN'(m_1,m_2,\ldots,m_n)\!=\!\cN'(m_2,\ldots,m_n)$]. In particular, 
$\cN'(2,\ldots,2)\!=\!\cN(2,\ldots,2)\!=\!1$.} 
\item{$\cK'(n,\mu)\!=\!\cK(n,\mu)\!\leq\!1 \Longleftrightarrow \mu\!\leq\!n+1$, 
and equality holds iff $\mu\!=\!n+1$. } 
\end{enumerate} 
\end{thm} 
\noindent
One should of course note the obvious fact that $\cN$ and $\cN'$ 
are symmetric functions in their arguments. Note also that conditions 
(a), (b), or (c) need not hold in assertions (2)--(4). 

\noindent 
{\bf Proof of Theorem \ref{thm:tri1}:} First note that the 
Newton polytopes must all be non-empty. 
The case (a) portion of assertions (0) and (1) then follows immediately 
from Corollary \ref{cor:zero}. Note also that the case (a) portion 
of assertion (0) implies the ``$\Longleftarrow$'' direction of 
assertion (2), since the underlying $n$-tuple of polytopes clearly has 
mixed volume zero. The ``$\Longrightarrow$'' direction of assertion (2) 
then follows easily from our earlier examples from Section \ref{sec:intro}.   
The case (b) portion of assertions (0) and (1) follows easily upon observing 
that $F$ is a linear system of $n$ equations in $n$ monomial terms, after 
multiplying the individual equations by suitable monomial terms. We can then 
finish by Proposition \ref{prop:easy}. 

To prove the case (c) portion of assertions (0) and (1), note that the case 
$n\!=\!1$ follows directly from UGDRS. For $n\!>\!1$, we have the following 
simple proof by induction: Assuming the desired bound holds for all 
$(n-1)\times (n-1)$ pyramidal systems, consider any $n\times n$ pyramidal 
system $F$. Then, via a 
suitable monomial change of variables, multiplying the individual equations 
by suitable monomials, and possibly reordering the $f_i$, we can assume that 
$f_1$ depends only on $x_1$. (Otherwise, $F$ wouldn't be pyramidal.) We thus 
obtain by UGDRS that $f_1$ has at most $m_1-1$ positive roots. By 
back-substituting these roots into $F'\!:=\!(f_2,\ldots,f_n)$, we obtain a 
new $(n'-1)\times (n'-1)$ pyramidal fewnomial system of type 
$(m'_2,\ldots,m'_{n'})$ with $n'\!\leq\!n$ and 
$m'_2\!\leq\!m_2,\ldots,m'_{n'}\!\leq\!m_{n'}$. 
By our induction hypothesis, we obtain that each such specialized $F'$ has at 
most $\prod^{n'}_{i=2}(m'_i-1)$ isolated roots in the positive orthant, and 
thus $F$ has at most $\prod^n_{i=1}(m_i-1)$ isolated roots in the positive 
orthant. (We already saw in the introduction that this bound 
can indeed be attained.) 

Our recursive formulae for $\cN'(2,m_2,\ldots,m_n)$ and 
$\cN(2,m_2,\ldots,m_n)$ from assertion (3) then follow by applying just the 
first step of the preceding induction argument, and noting that Proposition 
\ref{prop:easy} tells us that our change of variables preserves non-degenerate 
roots. 

Assertion (4) follows immediately from cases (a) and (b) of assertion (0). 
\qed 

One can of course combine and interweave families (a), (b), and (c) 
to obtain less trivial examples where we have exact formulae for 
$\cN(m_1,\ldots,m_n)$ and $\cK(n,\mu)$. More generally, 
one can certainly combine theorems \ref{thm:tri3} and \ref{thm:tri1} 
to obtain bounds significantly sharper than Khovanski's Theorem on 
Real Fewnomials, free from Jacobian assumptions, for many additional 
families of fewnomial systems. 

\section{Substitutions and Calculus: Proving Theorem 
\ref{thm:tri3} Minus Part (c)} 
\label{sec:tx} 
Let us preface our first main proof with some useful basic results. 
\begin{lemma}
\label{lemma:tri2} 
Let $F\!=\!(f_1,\ldots,f_n)$ be any $n\times n$ fewnomial 
system of type $(m_1,\ldots,m_n)$ with $m_1\!=\!1+\dim \newt(f_1)$. 
Then there is another $n\times n$ fewnomial system $G\!=\!(g_1,\ldots,g_n)$, 
also of type $(m_1,\ldots,m_n)$, such that  
$G$ has the same number of non-degenerate (resp.\ isolated) roots in 
$\Rn_+$ as $F$, $g_1:= 1\pm x_1\pm \cdots \pm x_{m_1-1}$ (with the signs in 
$g_1$ {\bf not} all ``$+$'') and, for all $i$, $g_i$ has $1$ as one of its 
monomial terms.  In particular, for $m_1\!=\!3$, we can assume further that 
$g_1\!:=\!1-x_1-x_2$. 
\end{lemma} 
 
\noindent 
{\bf Proof:} By dividing each $f_i$ by a suitable 
monomial term, we can assume that all the $f_i$ possess the monomial term $1$. 
In particular, we can also assume that the origin $\bO$ is a vertex of 
$\newt(f_1)$. Note also that the sign condition on $g_1$ must obviously hold, 
for otherwise the value of $g_1$ would be positive on the positive 
orthant. (The refinement for $m\!=\!3$ then follows by picking the monomial 
term one divides $f_1$ by more carefully.) So we now need only check 
that the desired canonical form for $g_1$ can be attained. 

Suppose $f_1\!:=\!1+c_1x^{a_1}+\cdots +c_{m_1-1}x^{a_{m_1-1}}$. 
By assumption, $\newt(f_1)$ is an $m_1$-simplex with vertex set 
$\{\bO,a_1,\ldots,a_{m_1-1}\}$, so $a_1,\ldots,a_{m_1-1}$ 
are linearly independent. Now pick any $a_{m_1},\ldots,a_n\!\in\!\Rn$ 
so that $a_1,\ldots,a_n$ are linearly 
independent. The substitution $x\mapsto x^{A^{-1}}$ (with $A$ the $n\times n$ 
matrix whose columns are 
$a_1,\ldots,a_n$) then clearly sends $f_1 \mapsto 1+c_1 x_1+\cdots 
+c_{m_1-1} x_{m_1-1}$, and Proposition \ref{prop:easy} tells us that 
this change of variables preserves degenerate and non-degenerate roots  
in the positive orthant. Then, via the change of variables 
$(x_1,\ldots,x_{m_1-1}) \mapsto (x_1/|c_1|,\ldots,x_{m_1-1}/|c_{m_1-1}|)$, we 
obtain that $g_1$ can indeed be chosen as specified. (The latter change 
of variables preserves degenerate and non-degenerate roots in the positive 
orthant for even more obvious reasons.) \qed 

Recall that a polynomial $p\!\in\!\R[x_1,\ldots,x_n]$ is {\bf homogeneous of 
degree $\pmb{D}$} iff 
$p(ax_1,\ldots,ax_n)\!=\!a^Dp(x_1,\ldots,x_n)$ for all $a\!\in\!\R$. 
\begin{prop} 
\label{prop:ind} 
Suppose $p\!\in\!\R[S_1,\ldots,S_n]$ is homogeneous of 
degree $D\!\geq\!0$. Also let $\alpha_1,u_1,v_1,\ldots,\alpha_n,
u_n,v_n\!\in\!\R$.
Then there is a homogeneous $q\!\in\!\R[S_1,\ldots,S_n]$, either identically 
zero or of degree \mbox{$D+n-1$,} such that 
$\frac{d}{dt}\left(p(u_1+v_1t,\ldots,u_n+v_nt)
\prod^n_{j=1}(u_j+v_jt)^{\alpha_j}\right)\!=\!q(u_1+v_1t,\ldots,u_n+v_nt)
\prod^n_{j=1}(u_j+v_jt)^{\alpha_j-1}$. 
\end{prop} 

\noindent 
{\bf Proof:} By the chain-rule, $\frac{d}{dt}\left(
p(u_1+v_1t,\ldots,u_n+v_nt) \prod^n_{j=1}(u_j+v_jt)^{\alpha_j}
\right)$ is simply\\ 
\scalebox{.9}[1]{$\left(\sum^n_{j=1}v_jp_j(u_1+v_1t,\ldots,u_n+v_nt)\right)
\left(\prod^n_{j=1}(u_j+v_jt)^{\alpha_j}\right)
+p(u_1+v_1t,\ldots,u_n+v_nt)\left(\sum^n_{i=1}\alpha_i v_i\frac{
\prod^n_{j=1} (u_j+v_jt)^{\alpha_j}}{u_i+v_it}\right)$}\\ 
where $p_i$ denotes the partial derivative 
of $p$ with respect to $S_i$. Factoring out a multiple of 
\mbox{$\prod^n_{j=1}(u_j+v_jt)^{\alpha_j-1}$} from the preceding expression, 
we then easily obtain that we can in fact take\\ 
$q(S_1,\ldots,S_n)\!=\!(v_1p_1(S_1,\ldots,S_n)+\cdots+
v_np_n(S_1,\ldots,S_n)) (S_1\cdots S_n) +p(S_1,\ldots,S_n)
\left(\sum^n_{i=1}\alpha_iv_i\frac{ S_1\cdots S_n}{S_i}\right)$. 
So we are done. \qed 

\begin{rolle} (1691) 
Let $g : [a,b] \longrightarrow \R$ be any continuous function with a 
derivative $g'$ well-defined on $(a,b)$. Then $g$ has $r$ roots in 
$[a,b] \Longrightarrow g'$ has at least $r-1$ roots in $(a,b)$. \qed 
\end{rolle}
\begin{lemma}
\label{lemma:unipert} 
Let $m\!\geq\!2$. 
Then for any real $c_1, u_1, v_1, \ldots,c_m,u_m,v_m$ and $[a_{ij}]$, the 
function
\[f(t):=\sum^m_{i=1}c_i\prod^n_{j=1}(u_j+v_jt)^{a_{ij}} \]
has no more than $n+\cdots+n^{m-1}$ roots in the open interval 
$I\!:=\!\{t\!\in\!\R_+ \; | \; u_j+v_jt\!>\!0 \text{ for all } j\}$. 

Furthermore, for any $\alpha_1,\ldots,\alpha_n\!\in\!\R$, $f$ has 
exactly $r$ roots in $I$ implies that there exist $\tilde{c}_1,\ldots,
\tilde{c}_m\!\in\!\R$ such that\\ 
\mbox{}\hfill 
$\tilde{f}(t)\!:=\!\sum^m_{i=1}\tilde{c}_i\prod^n_{j=1}(u_j+v_jt)^{a_{ij}}
\text{ has at least } r \text{ roots in } I, \text { no root of }
\tilde{f} \text{ in } I \text { is degenerate, and }$\\
\mbox{}\hfill$\text{ no root of } \tilde{f} 
\text{ in } I \text{ is an isolated root of } 
\left(\left(\prod^m_{j=1}(u_j+v_jt)^{\alpha_j}
\right)\tilde{f}'\right)'$.\hfill\mbox{} 
\label{lemma:ind}
\end{lemma}

\noindent 
{\bf Proof:} 
Throughout this proof let us consider only those roots lying in the open 
interval $I$ and assume that $f$ has exactly $r$ roots in $I$. We will 
in fact prove a stronger statement involving an extra parameter $D$ and then 
derive our lemma as the special case $D\!=\!0$.

First note that if 
\[g(t):=\sum^m_{i=1}p_i(u_1+v_1t,\ldots,u_n+v_nt)
\prod^n_{j=1}(u_j+v_jt)^{a_{ij}} \] 
for some homogeneous polynomials $p_1,\ldots,p_m$ of degree $D$, 
then 
\[g_0(t):=p_1(u_1+v_1t,\ldots,u_n+v_nt)+\sum^m_{i=2}
p_i(u_1+v_1t,\ldots,u_n+v_nt) \prod^n_{j=1}(u_j+v_jt)^{a_{ij}-a_{1j}} \] 
has the same number of roots in $I$ as $g$. 
In particular, using $D+1$ applications of Rolle's Theorem and Proposition 
\ref{prop:ind}, it is clear that $g_1\!:=\!g^{(D+1)}_0$ 
has at least $r-(D+1)$ roots, and we can in fact write 
\[g_1(t):=\sum^{m-1}_{i=1}q_i(u_1+v_1t,\ldots,u_n+v_nt)
\prod^n_{j=1}(u_j+v_jt)^{a'_{ij}}, \]
for some array $[a'_{ij}]$ and homogeneous polynomials $q_1,\ldots,q_m$ of 
degree $D+(D+1)(n-1)$.  

Now let $A(m,D)$ denote the maximum number of isolated roots of 
$g$ in the interval $I$. By what we've just observed, we immediately 
obtain the inequality 
\[ A(m,D)\leq A(m-1,nD+n-1)+D+1, \]  
valid for all $m\!\geq\!2$, $n\!\geq\!1$, and $D\!\geq\!0$. 
That $A(1,D)\!\leq\!D$ is clear, so one can then begin to 
bound $A(m,D)$ for general $m$ by recursion. A simple guess followed by an  
easy proof by induction yields 
\[ A(m,D)\leq (1+n+\cdots+n^m)(D+1)-1, \] 
which is valid for all $m,n\!\geq\!1$ and $D\!\geq\!0$. 
So the first assertion is proved. 

To prove the second part, note that the first part of our lemma implies 
that $f$ has only finitely many {\bf critical values} (i.e., values $f(x)$ 
with $f'(x)\!=\!0$) --- no more than $n+\cdots+
n^{m-1}$, in fact. Similarly, for any $\alpha_1,\ldots,\alpha_m\!\in\!\R$, 
there will only be finitely many roots for $\left(\left(
\prod^m_{j=1}(u_j+v_jt)^{\alpha_j} \right) f'\right)'$, unless 
this function is identically zero. In the latter case, no root of 
$\left(\left(\prod^m_{j=1}(u_j+v_jt)^{\alpha_j} \right) f'\right)'$ is 
isolated. So let us pick $\alpha_1,\ldots,\alpha_m\!\in\!\R$ so that 
$\left(\left(\prod^m_{j=1}(u_j+v_jt)^{\alpha_j} \right)f'\right)'$ is not 
identically zero. 

Note then that for all $\delta\!\in\!\Rs$ with $|\delta|$ 
sufficiently small, $f-\delta\prod^n_{j=1}(u_j+v_jt)^{a_{1j}}$ will have {\bf 
no} degenerate roots in $I$ and {\bf no} roots in $I$ making 
$\left(\left(\prod^m_{j=1}(u_j+v_jt)^{\alpha_j} 
\right) f'\right)'$ vanish. We can in fact guarantee that 
$f-\delta\prod^n_{j=1}(u_j+v_jt)^{a_{1j}}$ will also 
have at least $r$ non-degenerate roots in $I$ as follows: Let $n_+$ 
(resp.\ $n_-$) be the number of roots $t\!\in\!I$ of $f$ with $f'(t)\!=\!0$ 
and $f''(t)\!>\!0$ (resp.\ $f''(t)\!<\!0$). Clearly then, 
for all $\delta\!\in\!\Rs$ with $|\delta|$ sufficiently small, 
$f-\delta\prod^n_{j=1}(u_j+v_jt)^{a_{1j}}$ will have exactly $r+n_--n_+$ or 
$r+n_+-n_-$ roots in $I$, according as $\delta\!>\!0$ or $\delta\!<\!0$. 
(This follows easily upon dividing through by 
$\prod^n_{j=1}(u_j+v_jt)^{a_{1j}}$.) 
So let $\tilde{\delta}$ be sufficiently small, and of the correct sign, so 
that $f-\tilde{\delta}\prod^n_{j=1}(u_j+v_jt)^{a_{1j}}$ has at least $r$ roots 
in $I$, {\bf no} degenerate roots, and {\bf no} roots making 
$\left(\left(\prod^m_{j=1} (u_j+v_jt)^{\alpha_j} \right) f'\right)'$ vanish.

To conclude, simply let $\tilde{c}_1\!=\!c_1-\delta$ and 
$\tilde{c}_i\!:=\!c_i$ for all $i\!\geq\!2$. \qed 

\medskip 

\noindent 
{\bf Proof of Theorem \ref{thm:tri3} (Minus Part (c)):} 
We will reduce part (a) to part (b), prove part (b), 
and then refine our argument until we obtain 
$\cN'(3,3)\!=\!\cN(3,3)\!=\!5$. 

First note that in part (a), a simple Jacobian calculation 
reveals that the only way that $Z_+(f_1)$ can be degenerate is if $f_1$ is the 
square of a binomial. (Indeed, if $\newt(f_1)$ is a triangle then 
Lemma \ref{lemma:tri2} implies that $Z_+(f_1)$ is diffeomorphic to a 
line.) Part (0) of Theorem \ref{thm:tri1} then shows that 
our bound from (a) is easily satisfied in the special case where 
$f_1$ is a trinomial with $\newt(f_1)$ a line segment, so we 
can assume $\newt(f_1)$ is a triangle. Since $2+4+\cdots+2^{m-1}\!=\!2^m-2$, 
it then clearly suffices to prove part (b). 

To prove part (b), first note that UGDRS implies the case $n\!=\!1$, so 
we can assume $n\!\geq\!2$. Also, from the last paragraph, we already know 
that we can assume $\vol(\newt(f_1))\!>\!0$ when $n\!=\!2$. Since 
$F$ has no isolated roots when $n\!>\!2$ and the mixed volume 
of $\newt(f_1),\ldots,\newt(f_{n-1})$ is zero (via part (0) of theorem 
\ref{thm:tri1} again), we can assume henceforth that the mixed 
volume of $\newt(f_1),\ldots,\newt(f_{n-1})$ is positive. Since the supports 
of $f_1,\ldots,f_{n-1}$ can then all be translated into the vertex 
set of an $n$-simplex, Proposition \ref{prop:easy} tells us that we can 
assume in addition that $f_1,\ldots,f_{n-1}$ are {\bf affine} functions of 
$x_1,\ldots,x_n$.  
Letting $f_n(x_1,\ldots,x_n)\!=\!\sum^m_{i=1}c_i \prod^n_{j=1}x^{a_{ij}}_j$, 
we can then simply solve for $x_2,\ldots,x_n$ as functions of $x_1$ by 
applying Gaussian elimination to the first $n-1$ equations. Substituting into 
the last equation we then obtain a bijection between the roots of $F$ in the 
positive orthant and the roots of 
\mbox{$f(t):=\sum^m_{i=1}c_i \prod^n_{j=1}(u_j+v_jt)^{a_{ij}}
$} in the interval $I\!:=\!\{t\!\in\!\R_+ \; | \; u_j+v_jt\!>\!0 \text{ 
for all } j\}$, where $u_1,v_1,\ldots,u_n,v_n$ are suitable real constants. 

A simple Jacobian calculation then yields that $(\zeta_1,\ldots,\zeta_n)$ is 
a degenerate root of $F$ iff\\ 
\mbox{}\hfill $\left[ \sum^n_{\ell=1}v_\ell
\left. \frac{\partial f}{\partial x_\ell}\right|_{(\zeta_2,
\ldots,\zeta_n)=(u_2+v_2\zeta_1,\ldots,u_n+v_n\zeta_1)} 
\!=\!0 \text{ and } f(\zeta_1)\!=\!0\right]$,\hfill\mbox{}\\ 
and the above assertion is clearly true iff $f'(\zeta_1)\!=\!f(\zeta_1)\!=\!0$. 
So degenerate (resp.\ non-degenerate) roots of our univariate reduction 
correspond bijectively to degenerate (resp.\ non-degenerate) roots of $F$.  
Part (b) then follows immediately from Lemma \ref{lemma:ind}. 

To now prove that $\cN(3,3)\!=\!5$, thanks to Haas' 
counter-example, it suffices to show that $\cN(3,3)\!<\!6$. 
To do this, let us specialize our preceding notation to 
$(m,n)\!=\!(3,2)$, \ $(c_1,c_2)\!=\!(-A,-B)$,  
$(u_1,v_1,u_2,v_2)\!=\!(0,1,1,-1)$, and 
$(a_{11},a_{12},a_{21},a_{22})\!=\!(a,b,c,d)$, 
for some $a,b,c,d\!\in\!\R$ and positive $A$ and $B$. (Restricting 
$A,B,u_1,v_1,u_2,v_2$ as specified can easily be done simply by dividing $f_2$ 
by a suitable monomial term, as in the proof of Lemma \ref{lemma:tri2}.) 
In particular, the open interval $I$ becomes $(0,1)$. 

By using symmetry we can then clearly reduce to the following cases:\\ 
\begin{tabular}{ll}
{\bf A.} $a,b,c>0$ and $d<0$ \hspace{2cm}\mbox{} & {\bf B.} $a,c>0$ and 
$b,d<0$\hspace{1.85cm}\mbox{}\\
{\bf C.} $a,b>0$ and $c,d<0$ \hspace{2cm}\mbox{} & {\bf D.} $a,b,c,d>0$\hspace{3.1cm}\mbox{}\\
{\bf E.} $a,b,c,d<0$ \hspace{2cm}\mbox{} & {\bf F.} $a>0$ and $b,c,d<0$ \\ 
{\bf G.} $a,d>0$, $b,c<0$ & {\bf H.} At least one of the numbers $a,b,c,d$ 
is zero. 
\end{tabular}

Let $g(t)\!:=\!\frac{1}{B}t^{1-c}(1-t)^{1-d}f'(t)$. Then 
Lemma \ref{lemma:unipert} and our earlier substitution trick tells us 
that it suffices to show that any \[f(t):=1-A t^a(1-t)^b-Bt^c (1-t)^d, \] 
{\bf with all roots non-degenerate and no root of $f$ an  
isolated root of $g'$}, always has strictly less than $6$ roots 
in the open interval $(0,1)$. So let $r$ be the maximum number of roots in 
$(0,1)$ of any such $f$. 

Let us now prove $r\!<\!6$ in all 8 cases: 
\begin{itemize}
\item[{\bf A.}] $a,b,c>0$, $d<0$:  

Let $Q(x)=1-Ax^a(1-x)^b$ and $R(x)=Bx^c(1-x)^d$. 
The roots of $f$ may be regarded as the intersections in 
the positive quadrant of the parametrized curves $y\!=\!Q(x)$ and 
$y\!=\!R(x)$.  Since $\lim_{x\rightarrow 0^+} Q(x)=1$, 
$\lim_{x\rightarrow 1^-} Q(x)= 1$,   $\lim_{x\rightarrow 0^+} R(x)=0$, and 
$\lim_{x\rightarrow 1^-} R(x)=\infty$, it is easy to see via 
the Intermediate Value Theorem of calculus that the number of 
intersections must be odd. (One need only note that $f\!=\!Q-R$ and 
that the signs of $f'$ at the ordered roots of $f$ are nonzero and 
alternate.) So $r\!<\!6$.  
\item[{\bf B.}] $a,c>0$, $b,d<0$: \\
Almost exactly the same argument as case A will work here. The 
only difference here is that $\lim_{x\rightarrow 1^-} Q(x)\!=\!-\infty$.  
\item[{\bf C.}] $a,b>0$, $c,d<0$: \\See Lemma \ref{l1} below.
\item[{\bf D.}] $a,b,c,d>0$: \\See Lemma \ref{l2} below.
\item[{\bf E.}] $a,b,c,d<0$. \\ Multiplying $f(t)$ by 
$t^{\max\{-a,-c\}}(1-t)^{\max\{-b,-d\}}$, we can immediately reduce to 
case D.  
\item[{\bf F.}] $a>0$, $b,c,d<0$:\\See Lemma \ref{l3} below.
\item[{\bf G.}] $a,d>0$, $b,c<0$:\\See Lemma \ref{l4} below.
\item[{\bf H.}] At least one of the numbers $a,b,c,d$ is zero:\\
Use Lemma \ref{lemma:basic} below, noting that our hypotheses 
here imply that either $F$ or $\hat{F}$ is a 
quadratic polynomial. 
\end{itemize}
This concludes the proof of Theorem \ref{thm:tri3}, except for 
part (c), which we will complete in Section \ref{sec:tri}. \qed

\begin{rem} 
\label{rem:proj} 
Note that while we can attempt the same substitution trick 
for more complicated $F$, the complexity of the resulting 
recursion (involving derivatives and Rolle's Theorem) increases 
substantially. For instance, applying our proof in 
the special case where $n\!=\!2$ and $f_1(x,y)\!=\!1+x+cx^a-y$ unfortunately 
results in taking a number of derivatives which depends on $a$, thus 
obstructing a bound on the number of roots which is independent of the 
exponent vectors. \dia
\end{rem} 

We now detail the lemmata cited above. 
\begin{lemma}
\label{lemma:basic} 
Following the notation of the proof of Theorem \ref{thm:tri3}, recall 
that\\ $g(t)\!:=\!\frac{A}{B} t^{a-c}(1-t)^{b-d}(-a(1-t)+bt)-c(1-t)+dt$  
and that $r$ is the number of roots of
$f(t)\!:=\!1-At^a(1-t)^b-Bt^c(1-t)^d$ in 
the open interval $(0,1)$, where $f$ has no degenerate roots and 
no root of $f$ is an isolated root of $g'$. Also let\\ 
$F(u):=  -a(a-c)(a-c-1)u^3 
+(a-c) [2a(b-d+1)+b(a-c+1)]u^2+(d-b)[a(b-d+1)+2b(a-c+1)]u +b(b-d)(b-d-1)$, 
and\\ 
$\hat{F}(u):= -c(c-a)(c-a-1)u^3 +(c-a) [2c(d-b+1)+d(c-a+1)]u^2 
+(b-d)[c(d-b+1)+2d(c-a+1)]u +d(d-b)(d-b-1)$.  
Finally, let $N$ (resp.\ $M$) be the maximum number of non-degenerate roots in 
$(0,1)$ of $g$  (resp.\ the maximum of the number of positive roots of $F$ and  
$\hat{F}$), over all $(a,b,c,d)\!\in\!\R$ and $(A,B)\!\in\!\R^2_+$. 
Then $r-3\!\le\!N-2\!\le\!M\!\leq\!3$. 
\end{lemma} 

\noindent 
{\bf Proof:} Just as in the proof of Lemma \ref{lemma:ind},  
we easily see by Rolle's Theorem and division by suitable 
monomials in $t$ and $1-t$ that $r-1$ is no more than  
the number of roots in $(0,1)$ of $g$. So $r-1\!\leq\!N$. 
Note also that, in a similar way, $r-1$ is no more than the number 
of roots of $\hat{g}(t)\!:=\!\frac{B}{A}t^{c-a}(1-t)^{d-b}g(t)$ in 
$(0,1)$, and the latter function has the same number of roots (all of 
which are of course non-degenerate) in $(0,1)$ as $g$. 

To conclude, simply note that for suitable 
$\alpha,\beta,\gamma,\delta\!\in\!\R$, we have that 
$F\!\left(\frac{1-t}{t}\right)\!=\!t^\alpha(1-t)^\beta g''(t)$, 
$\hat{F}\!\left(\frac{1-t}{t}\right)\!=\!t^\gamma(1-t)^\delta \hat{g}''(t)$, 
and both expressions are cubic polynomials in $t$. 
So, by our preceding trick again, $N-2\!\leq\!M$, and thus 
$r-3\!\leq\!M$. That $M\!\leq\!3$ is clear from the fundamental 
theorem of algebra. \qed  

\begin{lemma}
Following the notation of Lemma \ref{lemma:basic}, let  
$T(x):=\frac{A}{B} x^{a-c}(1-x)^{b-d}(bx-a(1-x))$, $S(x):=c-(c+d)x$, 
$\hat{T}(x):=\frac{B}{A} x^{c-a}(1-x)^{d-b}(dx-c(1-x))$, and 
$\hat{S}(x):=a-(a+b)x$.  Then [$a,b>0$ and $c,d<0$] $\Longrightarrow r<6$. 
\label{l1}
\end{lemma}

\noindent 
{\bf Proof:} By Lemma \ref{lemma:basic}, we are done if $M\!<\!3$ or 
$N\!<\!5$. So let us assume $M\!=\!3$ to derive a contradiction. By Descartes' 
Rule of Signs (see Section \ref{sec:back} for a generalization), the 
coefficients of $F(u)$ or $\hat{F}(u)$ (ordered by exponent) must have 
alternating signs.  Thus, since $a,a-c,b,b-d\!>\!0$, we have that $a-c-1$ and 
$b-d-1$ must have the same sign. We then need to discuss two cases: 

\begin{itemize}
\item $a-c-1<0$ and $b-d-1<0$: 

This implies $c-a+1>0$ and $d-b+1>0$. Consequently, the 
coefficients of $u^3$ and $u^2$ in $\hat{F}(u)$ and $F(u)$ are all positive 
--- a contradiction.  
\item $a-c-1>0$ and $b-d-1>0$: 

The roots of $g$ in $(0,1)$ can be regarded as intersections of 
$y\!=\!T(x)$ and  $y\!=\!S(x)$, for $ 0\!<\!x\!<\!1$.
Since $T(\{0,1\})\!=\!0$, $-a(1-x)+bx\!=\!(a+b)x-a$ is strictly 
increasing, and $-a\!<\!0$, we must have that 
there is a smallest positive local minimum $c_0$ of $T$ with 
$T(c_0)\!<\!0$. Thus for $x$ near $c_0$, $T''(x)\!>\!0$. Since $T''(x)\!<\!0$ 
for $0\!<x\!\ll \!1$, there is $c^*\!\in\!(0,c_0)$ such that $T''(c^*)\!=\!0$.
Let $(x_1,y_1),\ldots,(x_K,y_K)$ be the intersection points
of $y\!=\!T(x)$ and $y\!=\!S(x)$
with $x_1<x_2<\cdots <x_K$. (A simple Jacobian calculation shows that 
$(x_i,y_i)$ is a degenerate root $\Longleftrightarrow x_i$ is a 
degenerate root of $g$. So every $(x_i,y_i)$ is in fact a non-degenerate 
root.) Then for all $i\!\in\!\{1,\ldots,K-1\}$ there is a 
$c_i\!\in\!(x_i,x_{i+1})$ with $T'(c_i)=-(c+d)>0$, and for all 
$i\!\in\!\{1,\ldots,K-2\}$ there is a $d_i\!\in\!(c_i,c_{i+1})$ with 
$T''(d_i)\!=\!0$.  Note that $c_0\!<\!c_1$. Thus $c^*\!<\!d_1$ and therefore 
$T''(x)\!=\!0$ has at least $K-1$ solutions. 
Since $T''$ and $F$ have the same number of 
positive roots (observing that $T''(u)/F(u)$ is a monomial in $u$ and 
$1-u$), we must have $N-1\!\le\!K-1\!\le\!3$. \qed 
\end{itemize}
\begin{lemma} Following the notation of Lemma \ref{l1},  
$a,b,c,d>0 \Longrightarrow r<6$. 
\label{l2}
\end{lemma}
{\bf Proof:} Again, by Lemma \ref{lemma:basic}, we need only show 
that $M\!<\!3$ or $N\!<\!5$. So let us assume $M\!=\!3$. Then by Descartes' 
Rule of Signs, $(a-c)(a-c-1)$ and $(b-d)(b-d-1)$ in the coefficients of $u^3$ 
and $u^0$ in $F(u)$ must have the same sign. There are now four cases to be 
examined.
\begin{itemize}
\item The signs of $a-c$, $a-c-1$, $b-d$, and $b-d-1$ are respectively 
$+,-,+$, and $-$: 

This makes the signs of coefficients of $u^3$ and $u^2$ of
$F(u)$ both positive. 

\item  The signs of $a-c$, $a-c-1$, $b-d$, and $b-d-1$ are respectively 
$-,-,+$, and $+$: 

Since  $b-d>0$, we have $d-b<0$ and $d-b-1<0$. This makes the constant term of
$\hat{F}(u)$ positive, and hence, the coefficients of $u$ and $u^2$ of 
$\hat{F}(u)$ must respectively be negative and positive. That is,
$c(d-b+1)+2d(c-a+1)<0\quad {\rm and} \quad 2c(d-b+1)+d(c-a+1)>0$. 
Thus, $-c(d-b+1)+d(c-a+1)<0.$ This is false, since $b-d-1>0$ and $a-c-1<0.$

\item The signs of $a-c$, $a-c-1$, $b-d$, and $b-d-1$ are all negative: 

By Descartes' rule of signs, $d-b-1$ and $c-a-1$ 
in the coefficients of $y^3$ and $y^0$ of $\hat{F}(y)$  
 must have the same sign. If both are negative, then coefficients of 
$u^3$ and $u^2$ of $F(u)$ would both be negative.
Thus  $d-b-1>0$ and $c-a-1>0$.
It is easy to see that $\hat{T}(x)\!<\!0$ for $0\!<x\!\ll\!1$ and
$\hat{T}(x)\!>\!0$ for $0\!<1-x\!\ll\!1$ and $\lim_{x\rightarrow 0^+} 
\hat{T}(x)=\lim_{x\rightarrow 1^-}\hat{T}(x)=0$.
Now let $L_0=\min \{c \; | \; 1\!>\!c\!>\!0, \ \hat{T}(c)\!<\!0 \text{ \ 
and \ } c \text{ \ is \ a \ local \ minimum of } \hat{T} \}$ and \\ 
$U_0=\max\{c \; | \; 1\!>\!c\!>\!L_0, \ c \text{ \ 
is \ a  \ local \ maximum of } \hat{T} \}$. 
Then for $x$ near $L_0$, $\hat{T}''(x)>0$. 
 Since $\hat{T}''(x)<0$ for $0<x\ll 1$, 
 there exists $L_1\!\in\!(0,L_0)$ such that $\hat{T}''(L_1)=0$.
Similarly, there is a $U_1\!\in\!(U_0,1)$ such that $\hat{T}''(U_1)\!=\!0$.

The roots of $\frac{B}{A}t^{c-a}(1-t)^{d-b}g$ can be studied via the 
intersections of $y=\hat{T}(x)$ and  $y=\hat{S}(x)$, for $ 0<x<1$.
Let $(x_1,y_1),\ldots,(x_k,y_k)$ be these intersection points, where  
$x_1\!<\!x_2\!<\cdots <\!x_k$. (A simple Jacobian calculation shows that 
$(x_i,y_i)$ is a degenerate root $\Longleftrightarrow x_i$ is a 
degenerate root of $g$. So every $(x_i,y_i)$ is in fact a non-degenerate 
root.) Then there are $c_1,\ldots,c_{k-1}$ with $c_i\!\in\!(x_i,x_{i+1})$ 
and $\hat{T}'(c_i)\!=\!-(a+b)<0$ for all $i\!\in\!\{1,\ldots,k-1\}$, 
and $d_1,\ldots,d_{k-2}$ with $d_i\!\in\!(c_i,c_{i+1})$ and 
$\hat{T}''(d_i)\!=\!0$ for all $i\!\in\!\{1,\ldots,k-2\}$.
If $x_1\!>\!L_0$, then $L_1<d_1$. If $x_1\!<\!L_0$, then $T(x_1)\!<\!0$. This 
implies $T(x_i)\!<\!0$ for all $i\!\in\!\{1,\ldots,k-2\}$, since the slope 
$-(a+b)$ of $\hat{S}(x)$ is negative. Therefore, $x_{k-2}\!<\!U_0$ and hence 
$d_{k-2}\!<\!U_1$.  So $\hat{T}''(x)\!=\!0$  has at least $k-1$ solutions. 
Since $\hat{T}''(x)\!=\!0$ and $\hat{F}(y)\!=\!0$ have the same number of 
solutions, we have $N-1\! \le \! k-1 \! \le \! M\!=\!3$. 

\item The signs of $a-c$, $a-c-1$, $b-d$, and $b-d-1$ are all positive: 
  
Since $a-c-1>0$ and $b-d-1>0$, the proof follows almost exactly the same line 
of reasoning as the last case, by intersecting the graphs of $T$ and 
$S$ instead of $\hat{T}$ and $\hat{S}$. \qed 
\end{itemize}

\begin{lemma} Following the notation of Lemma \ref{l1}, [$a>0$ and 
$b,c,d<0$] $\Longrightarrow r<6$.
\label{l3}
\end{lemma}
{\bf Proof:}  
Once again, by Lemma \ref{lemma:basic}, it suffices to show that 
$M\!<\!3$ or $N\!<\!5$. So let us assume that $M\!=\!3$. 
By checking the coefficients of $u^3$ and $u^0$ in $F(u)$, Descartes' Rule 
of Signs tells us that $a-c-1$ and $(b-d)(b-d-1)$ must have different signs. 
There are now three cases to be examined.
\begin{itemize}
\item  $a-c-1$, $b-d$, and $b-d-1$ are all negative: 
  
Then the signs of the coefficients of both $u^3$ and $u^2$ in  
$\hat{F}(u)$ will all be positive.
 
\item  The signs of $a-c-1$, $b-d$, and $b-d-1$ are respectively 
$-,+$, and $+$: 

Multiplying $f$ by $t^{-c}(1-t)^{-d}$ yields
$w(t)\!:=\!t^{-c}(1-t)^{-d}-At^{a-c}(1-t)^{b-d}-B$,  
where $-c\!>\!0$, $a-c\!>\!0$, $-d\!>\!1$, and $-d+b\!>\!1$. 
The roots of $w$ in $(0,1)$ can be regarded as the intersections of 
the parametrized curves $y\!=\!v(x)\!:=\!x^{-c}(1-x)^{-d}-Ax^{a-c}(1-x)^{b-d}$ 
and $y\!=\!B$. 
Let $(x_1,y_1),\ldots,(x_n,y_n)$ be the intersection points of
these two curves, where $x_1<x_2<\cdots <x_n$. 
(A simple Jacobian calculation shows that 
$(x_i,y_i)$ is a degenerate root $\Longleftrightarrow x_i$ is a 
degenerate root of $f$. So every $(x_i,y_i)$ is in fact a non-degenerate 
root.) Then for all $i\!\in\!\{1,\ldots,n-1\}$ there is a 
$c_i\!\in\!(x_i,x_{i+1})$ such that $v'(c_i)\!=\!0$. 
Thus $v'$ has at least $n-1$ roots in $(0,1)$. A straightforward 
computation then yields,\\
$v'(x):=Ax^{a-c-1}(1-x)^{b-d-1}(-(a-c)(1-x)+(b-d)x)
+x^{-c-1}(1-x)^{-d-1}(-c(1-x)+dx)$, 
which clearly has the same number of roots in $(0,1)$ as
\[t(x):=Ax^{a}(1-x)^{b}(-(a-c)(1-x)+(b-d)x)-c(1-x)+dx.\]
Thus $t''$ has at least $n-3$ roots in $(0,1)$. Since\\
$t''(x)/A=x^{a-2}(1-x)^{b-2}[-(a-c)a(a-1)(1-x)^3 
 +a((a+1)(b-d)+2(b+1)(a-c))x(1-x)^2$\\ 
\mbox{}\hspace{1.5cm}$-b((b+1)(a-c)+2(b-d)(a+1))x^2(1-x)+(b-d)b(b-1) x^3]$,\\ 
$t''$ has as many roots in $(0,1)$ as
\begin{equation*}
\begin{array}{ll}
P(u):=&-(a-c)a(a-1)u^3 +a((a+1)(b-d)+2(b+1)(a-c))u^2\\
    & -b((b+1)(a-c)+2(b-d)(a+1))u+(b-d)b(b-1) 
\end{array}
\end{equation*}
has positive roots. 
Since $a-1\!<\!a-c-1\!<\!0$, the coefficients of $u^3$ and $u^0$
in $P(u)$ are both positive. Thus $P$ has at most $2$ positive roots 
and we obtain $n-3\le 2$.

\item The signs of $a-c-1$, $b-d$, and $b-d-1$ are respectively $+,+$, and 
$-$: 
  
Since $a-c-1\!>\!0$ and $b-d\!>\!0$, it is easy to see that $T(x)\!<\!0$ for 
$0\!<\!x\!\ll\!1$ and $\lim_{x\rightarrow 1^-}T(x)\!=\!-\infty$. If $T(x)$ 
has no local minimum, then $y\!=\!T(x)$ and $y\!=\!S(x)$ have at most one 
intersection point. Otherwise, 
let $c_0\!:=\!\min \{c \; | \; 1\!>\!c\!>\!0, \ c \text{ \ is \ a \ local \ 
minimum \ of \ } T\}$. The rest of the proof is similar to that of Lemma 
\ref{l1}. \qed
\end{itemize}
\begin{lemma} Following the notation of Lemma \ref{l1}, 
[$a,d>0$ and $b,c<0$] $\Longrightarrow r<6$.
\label{l4}
\end{lemma}
{\bf Proof:}
One last time, Lemma \ref{lemma:basic} tells us that it suffices to 
prove that $M\!<\!3$ or $N\!<\!5$. So let's assume that $M\!=\!3$. Checking 
signs of coefficients of $u^3$ and 
$u^0$ of both $F(u)$ and $\hat{F}(u)$, Descartes' Rule of Signs tells us that 
$a-c-1<0$ and $d-b-1<0$. On the other hand, the alternating signs of 
coefficients of $u^2$ and $u^1$ of $F(u)$ yield 
$$2a(b-d+1)+b(a-c+1)<0 \quad {\rm and}\quad a(b-d+1)+2b(a-c+1)>0.$$ Thus,
$-a(b-d+1)+b(a-c+1)=a(d-1)+b(1-c)>0$. But this is impossible since 
$d-1\!<\!d-1-b\!<\!0$, $1-c\!>\!0$, $a\!>\!0$, and $b\!<\!0$. \qed

\begin{rem} 
When $A=1.12$, $B=0.71$, $a=0.5$, $b=0.02$, 
$c=-0.05$, and $d=1.8$, there are exactly $5$ roots of 
$1-Ax^a(1-x)^b-Bx^c(1-x)^d$ in $(0,1)$: They are, approximately,\\ 
\mbox{}\hfill\scalebox{1}[1]{$\{0.00396494, 0.02986317,0.4354707,0.72522344,
0.99620026\}$}.\hfill\mbox{}\\ 
In particular, this example is nothing more than the univariate reduction 
from the proof of Theorem \ref{thm:tri3} applied to a small 
perturbation of Haas' counter-example. \dia 
\end{rem}

\section{A Simple Geometric Approach, a Single Hard Case, and the 
Proof of Part (c) of Theorem \ref{thm:tri3} }
\label{sec:tri}
Let us begin with an extension of Rolle's Theorem to smooth curves
in the plane. 
\begin{lemma} 
\label{lemma:newrolle} 
Suppose $C$ is a smooth $1$-dimensional 
submanifold of $\R^2$ with:   
\begin{enumerate}
\item{At most $I$ inflection points that are isolated (relative  
to the locus of inflection points). } 
\item{At most $N$ non-compact connected components.} 
\item{At most $V$ isolated points of vertical tangency.} 
\end{enumerate}
Then the maximum finite number of intersections of any line with $C$ 
is $I+N+V+1$. 
\end{lemma} 

\noindent 
{\bf Proof:} Let $S^1$ be the realization of the circle obtained by 
identifying $0$ and $\pi$ in the closed interval $[0,\pi]$. 
Consider the natural map $\phi : C\longrightarrow S^1$ obtained by 
$x\mapsto \theta_x$ where $\theta_x$ is the angle in $[0,\pi)$ the normal line 
of $x$ forms with the $x_1$-axis. We claim that any $\theta\!\in\!S^1$ has at 
most $I+V+1$ pre-images under $\phi$. 

To see why, note that by assumption we can express $C$ as the union of 
no more than $I+V+1$ arcs where (a) any distinct pair of arcs is either 
disjoint or meets at $\leq\!2$ end-points, and (b) every end-point is either 
an isolated point of inflection or vertical tangency of $C$. 
(This follows easily by considering the graph whose vertices 
are the underlying inflection and vertical tangency points, 
and whose vertex adjacencies are determined by path-connectedness.) Calling 
these arcs {\bf basic arcs}, it is then clear that the interior of any basic 
arc is homeomorphic (via $\phi$) to a connected subset of 
$S^1\!\setminus\!\{0\}$. We then easily obtain that any 
$\theta\!\in\!S^1$ has at most $I+V+1$ pre-images under $\phi$, 
since each such pre-image belongs to exactly $1$ basic arc. 

Recall that a {\bf contact point} of a curve $C$ 
with a differential system 
$\frac{\partial \vec{X}}{\partial t}\!=\!\vec{G}(t)$  
is simply a point 
at which some solution of $\frac{\partial \vec{X}}{\partial t}\!=\!\vec{G}(t)$ 
has a tangent line in common with $C$. 
Now note that any line $L_m\!:=\!\{(x_1,x_n)\!\in\!\R^2 \; | \; 
m_1x_1+m_2x_2\!=\!m_0\}$ normal to $C$ forms an angle of  
$\mathrm{ArcTan}\!\left(\frac{m_2}{m_1}\right)$ with the $x_1$-axis. Thus, 
the number of contact points $C$ has with the differential system 
\[ \frac{\partial x_1}{\partial t}=m_2 \ , \ \frac{\partial x_2}
{\partial t}=-m_1 \] 
is at most $I+V+1$. By Rolle's Theorem for Dynamical Systems in the Plane 
(see, e.g., [Kho91,
corollary, pg.\ 23]), we then obtain that 
the number of intersections of $L_m$ 
with $C$ is at most $I+N+V+1$, 
for any real $(m_0,m_1,m_2)\!\neq\!(0,0,0)$. So we are done. \qed 

\begin{rem} 
The bound from Lemma \ref{lemma:newrolle} is tight in all cases. 
This is easily revealed by the examples in figure 
1 below and their obvious extensions. In particular, one can simply append 
$N-1$ disjoint lines to extend any example with $N\!=\!1$ to $N\!>\!1$. \dia   
\begin{figure}[h]
\epsfig{file=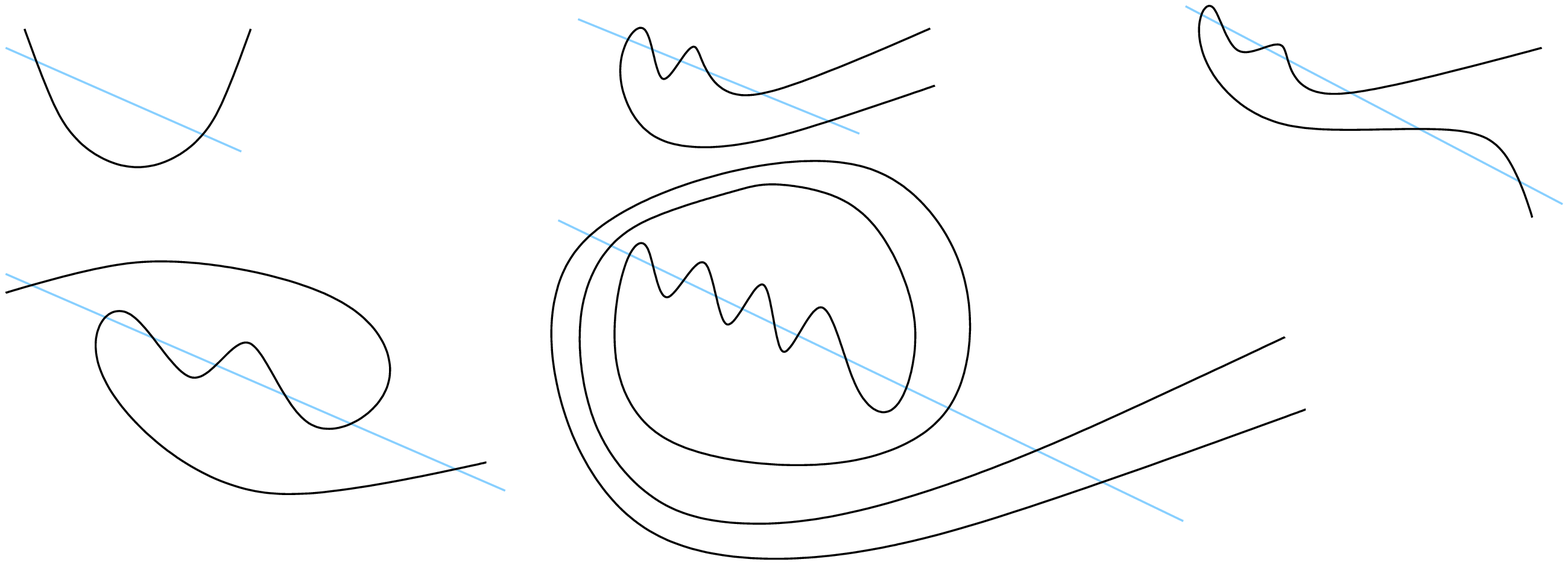,height=2.15in}
{\sc Figure 1} {\rm \scalebox{.95}[1]{
Lemma \ref{lemma:newrolle} gives a tight bound for $N\!=\!1$ 
and $(I,V)\!\in\!\{(0,0),(3,1),(4,1),(3,2),(7,5)\}$;}\\  
\mbox{}\hspace{1.75cm}and this generalizes easily to arbitrary $(I,N,V)$. } 
\end{figure} 
\end{rem}

We are now ready to give a concise geometrically motivated proof of 
the nearly optimal bound $\cN(3,3)\!\leq\!6$. 
This ``second'' proof of $\cN(3,3)\!\leq\!6$ was actually the original 
motivation behind this paper and, via a trivial modification,   
yields the proof of part (c) of Theorem \ref{thm:tri3} as well. 

\noindent
{\bf Short Geometric Proof of $\pmb{\cN(3,3)\!\leq\!6}$:}
Theorem \ref{thm:tri1} implies that we can assume that $f_1$ and $f_2$ have 
Newton polygons that are each triangles. Lemma \ref{lemma:tri2} of
the last section tells us that we can assume that $f_1\!=\!1\pm x_1\pm x_2$,
so we need only check the number of intersections of a line 
with $Z_+(f_2)$. In particular, since $Z_+(f_2)$ is diffeomorphic to a line 
(thanks to Proposition \ref{prop:easy}), Theorem \ref{thm:cool} tells us that 
$Z_+(f_2)$ has no more than $3$ inflection points and $1$ vertical tangent. 
By Lemma \ref{lemma:newrolle}, we are done. \qed  

\noindent
{\bf Proof of Part (c) of Theorem \ref{thm:tri3}:} 
Via a change of variables almost exactly like that of Lemma \ref{lemma:tri2} 
from the last section, we can assume that $f_1\!=\!1\pm x_1\pm x_2$. 
From here, we proceed exactly as in our 
last proof, noting that here $Z_+(f_2)$ instead has no more than 
$\area(\newt(p))$ isolated vertical tangents, $3\area(\newt(p))$ isolated 
inflection points, and $2D$ non-compact components (thanks to 
Theorem \ref{thm:cool}). \qed

Fewnomial curves happen to admit a simple ``fewnomial'' description
of their inflection points and singular points. This fact will be used here 
to prove our classification of when equality holds in our bound 
$\cN(3,3)\!\leq\!5$ (Corollary \ref{cor:poly}) and in the proof of Theorem 
\ref{thm:cool} in the next section. Let $\partial_i\!:=\!\frac{\partial}
{\partial x_i}$.
\begin{lemma}
\label{lemma:imp}
Suppose $f : \R^2_+ \longrightarrow \R$ is analytic.
Then [$z$ is an inflection point or a singular point of $Z_+(f)$]
$\Longrightarrow
\left\{f(z)\!=\!0 \text{ and } [\partial^2_1 f\cdot  (\partial_2 f)^2-
2\partial_1\partial_2 f\cdot\partial_1f\cdot\partial_2f
+\partial^2_2 f\cdot(\partial_1 f)^2]_{x=z}\!=\!0\right\}$. 
In particular, in the case where $f(x)\!:=\!p(x^{a_1},\ldots,x^{a_m})$
for some polynomial $p\!\in\!\R[S_1,\ldots,S_m]$
and $a_1,\ldots,a_m\!\in\!\R^2$,
the above cubic polynomial in derivatives is, up to a multiple which is a
monomial in $(x_1,x_2)$, a polynomial in $x^{a_1},\ldots,x^{a_m}$
with Newton polytope contained in $3\newt(p)$.
\end{lemma}
\noindent
{\bf Proof:} In the case of a singular point, the first assertion
is trivial. Assuming $\partial_2f\!\neq\!0$ at an inflection point then
a straightforward computation of $\partial^2_1x_2$ (via implicit
differentation and the chain rule) proves the first
assertion. If $\partial_2f\!=\!0$ at an
inflection point then we must have $\partial_1f\!\neq\!0$. So by computing
$\partial^2_2x_1$ instead, we arrive at the remaining case of the first
assertion. The second assertion follows routinely from the
chain rule. \qed

Let us now prove our polygonal classification of bivariate trinomial systems 
with maximally many roots in the positive quadrant.  

\noindent 
{\bf Proof of Corollary \ref{cor:poly}:} The segment case follows immediately 
from Corollary \ref{cor:zero}. For the remaining cases, Lemma 
\ref{lemma:tri2} implies that we can assume $f_1\!:=\!1\pm x_1 \pm x_2$ and 
$f_2\!:=\!1+Ax^a_1x^b_2+Bx^c_1x^d_2$ for some real $A$ and $B$. In particular, 
it is easily verified that the number of edges of $P_F$ and $P_G$ are the 
equal.  

So let $S_1\!:=\!Ax^a_1x^b_2$, $S_2\!:=\!Bx^c_1x^d_2$, and let 
$Z\!:=\!Z_+(f_2)$. Observe that Lemma \ref{lemma:imp} (along with a 
suitable rescaling of $f_2$ and the variables) tells us that we can 
bound the number of inflection points of $Z$ by analyzing the roots of a 
homogeneous polynomial in $(S_1,S_2)$ of degree $\leq\!3$. So let us now 
explicitly examine this polynomial in our polygonally defined cases. 

Clearly then, the triangle case corresponds to 
setting $a\!=\!d\!>\!0$ and $b\!=\!c\!=0$. 
We then obtain that [$x$ is an inflection point or a singular point of 
$Z$] $\Longrightarrow 1+S_1+S_2\!=\!0$ and $S_1+S_2\!=\!0$. So 
$Z$ has {\bf no} inflection points (or singularities). It is also even 
easier to see that $Z$ has no vertical tangents. So by Lemma 
\ref{lemma:newrolle}, 
$\cN(3,3)\!\leq\!2$ in this case. To see that equality can hold in this 
case, simply consider $F\!:=\!(x^2_1+x^2_2-25,x_1+x_2-7)$, 
which has $P_F\!=\!\conv(\{(0,0),(3,0),(0,3)\})$ and root set 
$\{(3,4),(4,3)\}$.  

For the quadrilateral case, we clearly have that $\newt(f_1)$ 
and $\newt(f_2)$ have exactly two inner edge normal vectors (with 
length $1$) in common. So let $v_i$ be the vertex of $\newt(f_i)$ 
incident to both the edges of $\newt(f_i)$ with these normals. 
Clearly then, we can assume that our above application 
of Proposition \ref{prop:easy} (which simply involved dividing 
the $f_i$ by suitable monomial terms and performing an invertible monomial 
change of variables) gives us $v_1\!=\!\bO$ as well. 
So we can assume $b\!=\!c\!=\!0$ and $a,d\!>\!0$. We then get the pair of 
equations $1+S_1+S_2\!=\!0$ and $a(d-1)S_1-d(a-1)S_2\!=\!0$,  
with $a,d\!\not\in\!\{0,1\}$. (If $\{a,d\}\cap\{0,1\}\!\neq\!\emptyset$ 
then $F$, or a suitable pair of linear combination of $F$, would be  
pyramidal and we would be done by Theorem \ref{thm:tri1}.)  
So $Z$ can have at most $1$ inflection point.  It is also even
easier to see that $Z$ has no vertical tangents. So by another application 
of Lemma \ref{lemma:newrolle}, $\cN(3,3)\!\leq\!4$ in this case. To see that 
equality can hold in this case, simply consider the system 
$(x^2_1-3x_1+2,x^2_2-3x_2+2)$, which has 
$P_F\!=\!\conv(\{\bO,(2,0),(2,2),(0,2)\})$ and root set 
$\{(1,1),(1,2),(2,1),(2,2)\}$. 

As for the pentagonal case, we can again assume (just as in the 
quadrilateral case) that our application of Proposition \ref{prop:easy} 
placed the correct vertex of $\newt(f_1)$ at the origin. In particular,  
we can assume $b\!=\!0$ and $a,c,d\!>\!0$. We then get 
the pair of equations $1+S_1+S_2\!=\!0$ and 
$a^2(d-1)S^2_1+a(ad-d-2c)S_1S_2-c(c+d)S^2_2\!=\!0$,   
with $ac(d-1)(c+d)\!\neq\!0$. (Similar to the last case, it is easily checked 
that if the last condition were violated, then we would be back in one of our 
earlier solved cases.) However, a simple check of the discriminant of the 
above quadratic form in $(S_1,S_2)$ shows that there is at most $1$ 
root, counting multiplicities, in any fixed quadrant. So, similar to the 
last case, we obtain $\cN(3,3)\!\leq\!4$ in this case. 
To see that the equality can hold in this case, simply 
consider the system $(x^2_2-7x_2+12,-1+x_1x_2-x^2_1)$, which has  
$P_F\!=\!\conv(\{\bO,(2,0),(2,2),(1,3),(0,2)\})$ and  
root set $\left\{\left(3,\frac{3\pm\sqrt{5}}{2}\right),\left(4,2\pm\sqrt{3}
\right)\right\}$. \qed  

\section{Monomial Morse Functions and Connected Components: Proving Theorem 
\ref{thm:cool}} 
\label{sec:morse}
Let us begin with a refinement of a lemma due to Khovanski.\footnote{In 
\cite{risler}, Risler outlines a proof of the first portion of Lemma 
\ref{lemma:kho}, in the special case where all the $a_i$ lie in $\Zn$. He 
also cites a paper of Khovanski for further details. However, as far as the 
authors can tell, the wrong paper by Khovanski was cited.} Recall that a {\bf 
$\pmb{d}$-flat} in $\Rn$ is simply a translate of a $d$-dimensional subspace 
of $\Rn$. 
\begin{lemma}
\label{lemma:kho} 
Let $a_1,\ldots,a_m\!\in\!\Rn$, $c_1,\ldots,c_m\!\in\!\R$, and let $Z$ be 
the zero set of $\sum^m_{i=1}c_ie^{a_i\cdot z}$ in $\Rn$, 
where $\cdot$ denotes the usual Euclidean inner product of vectors in $\Rn$. 
Then there is an $(n-1)$-flat in $\Rn$ which intersects at least half of the 
non-compact connected components of $Z$. Furthermore, the set of 
unit normal vectors of all such $(n-1)$-flats is non-empty and open in the 
unit $(n-1)$-sphere. 
\end{lemma} 
\noindent
For the convenience of the reader, we sketch a proof below, based 
on an argument of Jean-Jacques Risler from \cite[Sec.\ 2]{risler}. 

\noindent
{\bf Proof:} Let us call any function of the form described in 
the statement of the theorem an {\bf exponential $\pmb{m}$-sum}. 
Let $Z_i$ be any non-compact component of $Z$ and $C_i$ 
any connected unbounded curve (defined by a system of exponential 
$m$-sums) lying in $Z_i$. Let $p_i$ be any limit point (as 
$\|z\|\longrightarrow +\infty$) of the set 
$\left\{\left. \frac{z}{\|z\|} \; \right| \; z\!\in\!C_i\right\}$. That the 
set of such limit points is in fact finite follows easily from a slightly more 
general version of Khovanski's Theorem on Real Fewnomials \cite[Cor.\ 6, Pg.\ 
80, Sec.\ 3.12]{few}, stated in terms of exponential sums.   

If $H$ is any hyperplane (so $\bO\!\in\!H$) such that 
$p_i\!\not\in\!H$ for all $i$ then one of the open unit hemispheres 
defined by $H$ contains at least half of the points $p_i$. 
In particular, note that such an $H$ must clearly exist and 
any hyperplane $H'$ with unit normal vector sufficiently near that of
$H$ will also define an open unit hemisphere containing at least 
half the $p_i$. To conclude, note that any $(n-1)$-flat, parallel 
to $H$ and far enough in the direction of the pole of the hemisphere 
containg the most $p_i$, will intersect half of the $C_i$ and thus half of 
the $Z_i$.  So we are done. \qed 

We will also need an extension of the classical bounds on the 
number of connected components of a real algebraic set. 
\begin{lemma} 
\label{lemma:finite} 
Given any $\mu$-sparse $k\times n$ fewnomial system $F$, 
the number of connected components of $Z_+(F)$ is 
no more than $2^{n-\frac{1}{2}}(2n+1)^\mu 2^{\mu(\mu+1)/2}$. \qed 
\end{lemma} 

\noindent 
The smooth case, which admits a sharper bound, is detailed in 
\cite[Sec.\ 3.14]{few}. The special case of integral exponents 
(allowing degeneracy) is nothing more than \cite[Cor.\ 3.2]{real} and 
the proof in \cite[Sec.\ 3.2]{real} extends with
no difficulty to real exponents. One can in fact generalize the 
above lemma to semi-Pfaffian sets, provided one loosens the stated 
upper bound somewhat \cite{zell}. 

A construction which will prove quite useful when we count 
components via critical points of maps is to find a monomial map which is a  
{\bf Morse function} \cite{morse} relative to a given fewnomial zero set. 
Recall that an {\bf $\pmb{n}$-dimensional polyhedral cone}  
in $\Rn$ is simply a set of the form $\{r_1a_1+\cdots+r_na_n \; | \; 
r_1,\ldots, r_n\!\geq\!0\}$, 
where $\{a_1,\ldots,a_n\}$ is a generating set for $\Rn$. In particular, 
an $n$-dimensional cone in $\Rn$ always has non-empty interior. 
\begin{lemma} 
\label{lemma:fiber}
Suppose $f$ is an $n$-variate $m$-nomial. 
Then there is an $n$-dimensional polyhedral cone
$K\!\subseteq\!\Rn$ such that $a\!\in\!K\!\setminus\!\{\bO\}$ implies 
\begin{enumerate}
\item{Every critical point of the restriction of $x^a$ to $Z_+(f)$ 
is non-degenerate. } 
\item{The level set in $Z_+(f)$ of any regular value of $x^a$ has 
dimension $\leq\!n-2$.} 
\item{No connected components of $Z_+(f)$ other than 
isolated points are contained in any level set of $x^a$.}
\item{At least half of the non-compact connected component of $Z_+(f)$ have 
unbounded values of $x^a$.}  
\end{enumerate} 
\end{lemma}

\noindent 
{\bf Proof:} Set $Z\!:=\!Z_+(f)$. Let us first show how assertion (3) 
can be attained: Since the 
number of components of $Z$ is finite by Lemma \ref{lemma:finite},
we can temporarily assume that $Z$ consists of a single component. 
Then, if we could find $n$ linearly independent $a\!\in\!\Rn$ with 
$Z\!\subset\!\{x\!\in\!\Rn_+ \; | \; x^a\!=\!c_a$\} for some $c_a$, 
Proposition \ref{prop:easy} would immediately imply that $Z$ is contained 
in a point. So condition (3) can be enforced.  

To ensure the truth of condition (4), note that we can perform the 
substitution $(x_1,\ldots,x_n)\!=\!(e^{z_1},\ldots,e^{z_n})$ 
to reduce to finding an $(n-1)$-flat intersecting 
at least half of the non-compact components of the real zero 
set of an exponential $m$-sum. That such a hyperplane exists (and the fact 
that a small open neighborhood of such hyperplanes exist) then follows 
from Lemma \ref{lemma:kho}, so condition (4) holds. 

To enforce conditions (1) and (2) let us maintain our 
last change of variables. A simple derivative computation (noting 
that $x\mapsto (e^{x_1},\ldots,e^{x_n})$ defines a 
diffeomorphism from $\Rn$ to $\Rn_+$) then shows that 
it suffices to instead prove the analogous statement where $f$ is replaced by 
a real exponential sum (a real analytic function in any event) and $x^a$ is 
replaced by the linear form $a_1z_1+\cdots+a_nz_n$. The latter analogue is 
then nothing more than an application of \cite[Lemma 1, Pg.\ 304]{bcss}. 
Since the number of components of $Z$ is finite, we thus obtain 
assertions (1) and (2). 

Our lemma then follows by intersecting the four sets of $a$ 
we have just determined, and noting (thanks to Lemma \ref{lemma:kho}) 
that the intersection clearly contains a cone over a small $(n-1)$-simplex. 
\qed 

\noindent 
Note that if we omit condition (4) then we can make conditions (1)--(3) 
hold in an even larger set of $a$: all $a$ in $\Rn$ outside a finite union 
of hyperplanes. 

\begin{ex} 
In general, one can not find an $a$ with {\bf every} component of $Z_+(f)$ 
giving unbounded values for $x^a$. This follows from an elementary calculation 
with $n\!=\!2$ and 
\[f(x,y)\!:=\!\left(1-x-xy-\frac{1}{y}\right)\left(1-y-xy-\frac{1}{x}\right)
\left(1-\frac{1}{x}-\frac{1}{y}\right),\] showing that, for all $(a_1,a_2)$, 
$x^{a_1}y^{a_2}$ has unbounded values on no more than $2$ of the $3$ 
components of $Z_+(f)$. The authors thank Daniel Perrucci and Fernando Lopez 
Garcia for this example. \dia 
\end{ex} 

We will also need the following useful perturbation result, which 
can be derived via Sard's Theorem \cite{hirsch} and a simple homotopy 
argument. (See, e.g., 
\cite[Lemma 2]{basu} for even stronger results of this form in the setting of 
integral exponents and zero sets in $\Rn$.) 
\begin{lemma}
\label{lemma:pert} 
Let $f$ be any $n$-variate $m$-nomial, $Z^\delta_+(f)$ the 
solution set of $|f|\!\leq\!\delta$ in $\R^n_+$, and 
$\stackrel{\circ}{Z^\delta_+}(f)$ the boundary of $Z^\delta_+(f)$. 
Then for $\delta\!>\!0$ 
sufficiently small, $\stackrel{\circ}{Z^\delta_+}(f)$ is smooth and  
has at least as many compact (resp.\ non-compact) connected components as 
$Z_+(f)$. \qed 
\end{lemma}
\begin{ex} 
If $f(x,y)\!:=\!x^2+(1-xy)^2$ then note that 
$Z_+(f)$ is empty while, for any $\delta\!>\!0$, $Z^\delta_+(f)$ 
contains the point $\left(\sqrt{\delta},\frac{1}{\sqrt{\delta}}\right)$. So 
$Z^\delta_+(f)$ need not have the exactly same number 
of compact (or non-compact) components 
as $Z_+(f)$, even if $\delta\!>\!0$ is very small. 
The authors thank Daniel Perrucci for this example. \dia 
\end{ex}  

\begin{ex} 
Boundaries of tubes about analytic sets behave a bit 
differently in $\Rn_+$ than in $\Rn$. For instance, unlike 
the analogous bound over $\Rn$ (see, e.g., \cite[Lemma 3.1]{real}),   
the number of components of $Z_+(f)$ can {\bf not} be 
bounded above by {\bf half} the number of components of $\stackrel{\circ}
{Z^\delta_+}(f)$: taking $f(x,y)\!=\prod^D_{i=1}(y-ix)$, it is easily checked 
that $Z_+(f)$ has exactly $D$ components, while $\stackrel{\circ}
{Z^\delta_+}(f)$ has exactly $D+1$ components for $\delta\!>\!0$ sufficiently 
small. \dia 
\end{ex}

\noindent 
{\bf Proof of Theorem \ref{thm:cool}:} 
For simplicity, let us assume that all components are connected and lie in 
$\Rn_+$. 

Assertion (0) follows immediately 
from Proposition \ref{prop:easy}, UGDRS, and noting that $0$ is an $n$-variate 
polynomial with exactly one non-compact connected component in its real zero 
set. 

To prove assertion (1), note that assertion (2) of Theorem \ref{thm:cool} 
(which we'll soon prove below) and assertion (4) of Theorem \ref{thm:tri1} 
easily imply the first formula. Proposition \ref{prop:easy} then 
tells us that the second equality can be proved simply by employing a 
monomial change of variables to reduce to the case of an $(m-1)$-variate 
$m$-nomial $g$. In particular, since $m\!\geq\!3$, every component  
of $Z_+(g)$ will still be non-compact. (This is clear from 
another application of Proposition \ref{prop:easy}, separating 
the cases where $\dim \newt(g)$ equals, or is strictly less than, $m$.) 
Moreover, $\dim \newt(g)\!=\!m \Longrightarrow $ the number of 
non-compact components is exactly $1$. So we can assume that 
$\dim \newt(g)\!<\!m$, use Proposition \ref{prop:easy} one last time, 
and then intersect with an appropriate coordinate flat to derive 
the final inequality. So assertion (1) is proved and we can assume henceforth 
that $n\!\geq\!2$. 

To prove assertion (2), let us first construct a concrete family of examples 
realizing the lower bound: Consider the polynomials  
\[ g_1(x):= \left(\sum^n_{i=2}(x_i-1)^2\right)+
\left(\prod^{\lfloor m/2\rfloor -n-1}_{i=1}(x_1-i)^2\right) 
\ \ \ \  \text{ and }  \ \ \ \ \ 
g_2(x):= \sum^n_{j=1} \ \prod^{\lfloor (m-1)/(2n)\rfloor}_{i=1} 
(x_j-i)^2. \] 
Clearly, $g_1$ and $g_2$ respectively have exactly 
$2\lfloor m/2\rfloor$ and $2n\left\lfloor (m-1)/(2n)\right\rfloor$ monomial 
terms. From the basic fact that $a^2+b^2\!=\!0 \Longrightarrow 
a\!=\!b\!=\!0$ for all real $(a,b)$, it is easily checked that 
the numbers of roots of $g_1$ and $g_2$ in $\Rn_+$ are finite and 
in fact are identically the formulae embedded in our lower bound. More 
precisely, we immediately obtain 
our lower bound {\bf for a restricted class of $\pmb{m}$ depending 
on the congruence class of $\pmb{m}$ mod $\pmb{2}$ or mod $\pmb{2n}$.} 
This restriction can easily be removed by adding additional monomial 
terms in such a way that the number of compact and non-compact components
is not decreased. To do this, simply note that by Sard's Theorem 
\cite{hirsch} (and the definition of an $n$-sphere), we have that 
$Z_+(g_i-\delta_0)$ is smooth and has the same numbers of 
compact and non-compact components as $Z_+(g_i)$ for all 
$\delta_0\!>\!0$ sufficiently small. Similarly, the same will then be true 
of $Z_+(x_1(g_i-\delta_0)-\delta_1)$, for $\delta_1\!>\!0$ sufficiently 
small, and the latter polynomial has exactly $1$ more monomial term 
than $g_i$. Proceeding inductively, we can thus remove our restriction 
on $m$, and we thus obtain the lower bound of assertion (2). 

To prove the upper bound of assertion (2), note that we can 
divide by a suitable monomial so that $f$ has a nonzero constant 
term. By Lemma \ref{lemma:pert}, we then have that for $\delta\!>\!0$ 
sufficiently small, it suffices to bound the number of compact 
components of $\stackrel{\circ}{Z^\delta_+}(f)$ --- an ``envelope'' of 
$Z_+(f)$. Recall that Lemma \ref{lemma:pert} also grants us that 
$\stackrel{\circ}{Z^\delta_+}$ can be assumed to be smooth. 

By Proposition \ref{prop:easy} and Lemma \ref{lemma:fiber}, we can 
then pick an $n\times n$ matrix $A$ so that, after we 
make the change of variables $x\!=\!y^A$, the number of compact and 
non-compact components of 
$\stackrel{\circ}{Z^\delta_+}$ is preserved, {\bf no} component of 
$\stackrel{\circ}{Z^\delta_+}$ of positive dimension is contained 
in a hyper-plane parallel to the $y_1$-coordinate hyperplane, 
and we can use critical points of the function $y_1$ to count compact 
components.  

Consider then the systems of 
equations $G_\pm\!:=\!(f\pm\delta,y_2\partial_2f,\ldots,y_n\partial_nf)$,  
where $\partial_i\!=\!\frac{\partial}{\partial y_i}$ here. 
By construction, every compact component of $\stackrel{\circ}{Z^\delta_+}(f)$ 
results in at least two extrema of the function $y_1$, i.e., 
$P_{\mathrm{comp}}(n,m)$ is bounded above by an integer no greater than 
half the total number of roots of $G_+$ and $G_-$.  (In particular, if 
$Z_+(f)$ were smooth to begin with, then we could have omitted the use 
of $\stackrel{\circ}{Z^\delta_+}(f)$ and $G_\pm$, since 
$P_{\mathrm{comp}}(n,m)$ 
would instead be bounded above by an integer no greater than half the number 
isolated roots of $G\!:=\!(f,y_2\partial_2f,\ldots,y_n\partial_nf)$.) 
Note also that by assertion (1) of Lemma \ref{lemma:fiber}, all the roots of 
$G_\pm$ (and $G$) are non-degenerate. Furthermore, $G_\pm$ (and $G$) clearly 
has no more than $m$ distinct exponent vectors, so the upper bound 
on $P_{\mathrm{comp}}(n,m)$ holds. As for the number of compact 
components of $Z_+(\rho)$, the preceding argument applies as well, 
so we need only observe that $\rho\pm\delta$ and $y_2\partial_2\rho$ 
are both of the form $q(x^{r_1}y^{s_1},x^{u_2}y^{v_2})$ where 
$\newt(q)\!=\!\newt(p)$. So assertion (2) follows. 

To prove assertion (3), 
let us construct another family of explicit examples: Consider the 
polynomials  
\[ h_1(x):= \prod^{m-1}_{i=1}(x_1-i) \ \ \ \ \ \text{ and } \ \ \ \ \ 
h_2(x):= \sum^{n-1}_{j=1} \ \ \prod^{\lfloor (m-1)/(2n-2)\rfloor}_{i=1} 
(x_j-i)^2. \] 
Clearly, $h_1$ and $h_2$ respectively have exactly 
$m$ and $2(n-1)\lfloor (m-1)/(2n-2)\rfloor$ monomial 
terms. Note also that $Z_+(h_1)$ and $Z_+(h_2)$ 
have only non-compact connected components, and the numbers of such  
components are in fact the formulae embedded in our lower bound (for a 
restricted class of $m$). The lower bound of assertion (3) then follow easily, 
mimicking the argument we used earlier to remove the congruence 
class restriction which arose during the proof of the lower bound of assertion 
(2).  

To prove the upper bounds of assertion (3), 
let us work directly with $f$ and make an 
independent application of Lemma \ref{lemma:fiber}. We can then  
apply Proposition \ref{prop:easy} to make a change of variables $x\!=\!y^{A'}$ 
(preserving the number of compact and non-compact components of $Z_+(f)$), 
so that at least half of the non-compact components of
$Z_+(f)$ have unbounded values of $y_1$. So, for $\eps\!>\!0$ 
sufficiently small, the number of non-compact components of $Z_+(f)$ is no 
more  than twice the number of components 
of the intersection $Z'\!:=\!Z_+(f)\cap\!\left\{y\!\in\!\Rn_+ \; | 
\; y_1\!=\!\frac{1}{\eps}\right\}$. 
So by substituting $y_1\!=\!\frac{1}{\eps}$ into $f$, we obtain a new 
$m$-nomial hypersurface $Z''\!\subseteq\!\R^{n-1}$ 
with at least as many components as $Z'$. So $Z''$ has at least 
half as many components as $Z_+(f)$ has non-compact components, and thus the 
number of non-compact components of $Z_+(f)$ is no more than $2P(n-1,m)$. 
So the upper bound 
on $P_{\mathrm{non}}(n,m)$ is proved. As for the number of non-compact 
connected components of $Z_+(\rho)$, the preceding argument still applies. 
So we need only observe that, modulo a monomial change of variables via 
Proposition \ref{prop:easy}, $\rho$ can be assumed to be a polynomial 
of degree $D$. Lemma \ref{lemma:kho} and B\'ezout's Theorem 
\cite[ex.\ 1, pg.\ 198]{sha} (along with an exponential 
change of variables) then proves what is left of assertion (3). 
 
To prove assertion (4) simply note that the isolated points of vertical 
tangency of $Z_+(f)$ are exactly the isolated roots of the bivariate fewnomial 
system $H\!:=\!\left(f,x_2\frac{\partial f}{\partial x_2}\right)$. When 
$f(x)\!=\!p(x^{a_1},\ldots,x^{a_m})$ for some $p\!\in\!\R[S_1,
\ldots,S_m]$, a simple application of the 
chain rule then shows that $x_2\frac{\partial f}{\partial x_2}\!=\!
q(x^{a_1},\ldots,x^{a_m})$ for some $a_1,\ldots,a_m\!\in\!\R^2$ and some 
$q\!\in\!\R[S_1, \ldots,S_m]$ with $\newt(q)\!\subseteq\!
\newt(p)$. In particular, $p\!=\!1+S_1+\cdots +S_m \Longrightarrow 
q$ is a homogeneous linear form in $S_1,\ldots,S_m$, so 
the first part of assertion (4) follows. To prove the 
second part, note that $(r_1,s_1)$ and $(u_2,v_2)$ linearly 
dependent $\Longrightarrow Z_+(\rho)$ is a union of no more than two 
binomial curves (via Proposition \ref{prop:easy} and factoring over 
$\R$), and such curves have no isolated points of vertical tangency 
in $\R^2_+$. So, assuming $\det \begin{bmatrix}  
r_1 & u_2 \\ s_1 & v_2 \end{bmatrix}\!\neq\!0$, Proposition 
\ref{prop:easy} then tells us that it suffices to count the 
isolated roots $(S_1,S_2)\!\in\!\R^2_+$ of the $2\times 2$ polynomial system 
$H$. By Bernstein's Theorem \cite{bkk}, the number of complex 
isolated roots of the resulting system is at most $\area(\newt(p))$. 
So assertion (4) is proved. 

To prove assertion (5), first note that [$m\!=\!1 \Longrightarrow 
Z_+(f)$ is empty] and [$m\!=\!2 \Longrightarrow Z_+(f)$ has 
no isolated inflection points]. So we can assume $m\!=\!3$ and that 
$f$ has a constant term.  Note then that by Lemmata 
\ref{lemma:tri2} and \ref{lemma:imp}, $(x_1,x_2)$ is an inflection point of 
$Z_+(f)  \Longrightarrow f(x)\!=\!q(x^{r_1}y^{s_1},x^{u_2}y^{v_2})\!=\!0$, 
where $q\!\in\!\R[S_1,S_2]$. In particular, 
Lemma \ref{lemma:imp} tells us that $q$ is either a homogeneous cubic 
polynomial or a polynomial with Newton polytope contained in 
$3\newt(p)$, according as we focus on the first or second part of assertion 
(5). Just as in the last paragraph, we can also assume that 
$\det \begin{bmatrix} r_1 & u_2 \\ s_1 & v_2 \end{bmatrix}\!\neq\!0$, 
and thus reduce to counting the isolated roots in $\R^2_+$ of a $2\times 2$ 
polynomial system in $(S_1,S_2)$. For the first part of assertion (5), 
the fundamental theorem of algebra tells us that $q$ splits completely 
over $\C[S_1,S_2]$, so we can further reduce to no more than three 
$2\times 2$ linear systems and easily obtain our bound of $3\cK'(2,m)$. 
For the second part, we can easily conclude by Bernstein's Theorem \cite{bkk}. 
So assertion (5) is proved. 

To prove the last observation of Theorem \ref{thm:cool}, note that 
by Proposition \ref{prop:sandwich} and Lemma \ref{lemma:tri2}, 
$Z_+(f)$ has a singularity $\Longrightarrow \newt(f)$ is a line segment, 
and then $f$ must be the square of a binomial. So the case where $Z_+(f)$ is 
singular follows immediately. The case where $Z_+(f)$ is smooth then follows 
easily from assertions (4) and (5), since Theorem 
\ref{thm:tri1} of Section \ref{sec:back} implies that 
$\cK'(2,3)\!=\!\cK(2,3)\!=\!1$. \qed  

\section{Momenta, Polytopes, and the Proof of Theorem \ref{thm:moment} } 
\label{sec:moment} 
Let $\bar{S}$ and $\inte(S)$ respectively denote the topological closure 
and topological interior of any set $S$, and let 
$\relint(Q)$ denote the relative interior of any $d$-dimensional 
polytope $Q\!\subset\!\Rn$, i.e., $Q\!\setminus\!R$ where 
$R$ is the union of all faces of $Q$ of dimension strictly 
less than $d$ (using $\emptyset$ as the only face of 
dimension $<\!0$). We then have the following variant of the 
momentum map from symplectic geometry \cite{smale,souriau}. 
\begin{lemma}
\label{lemma:moment}
Given any $n$-dimensional convex compact polytope $P\!\subset\!\Rn$, there 
is a real analytic diffeomorphism $\psi_P : \Rn_+ \longrightarrow 
\inte(P)$. In particular, if $f$ is an $n$-variate $m$-nomial with 
$\newt(f)\!=\!P$ (so $\dim \newt(f)\!=\!n$) and 
$w\!\in\!\Rn\!\setminus\!\{\bO\}$, then 
$\psi_P(Z_+(f))$ has a limit point in $\relint(P^w) \Longrightarrow 
\init_w(f)$ has a root in $\Rn_+$. Moreover, there is a 
real analytic diffeomorphism between $Z_+(\init_w(f))\!\subset\!\Rn_+$ and 
$\left(\relint(P^w)\cap\overline{\psi_P(Z_+(f))}\right)
\times \R^{n-\dim P^w}_+$. 
\end{lemma} 

\noindent
{\bf Proof:} By \cite[Sec.\ 4.2, Lemma, Pg.\ 82]{tfulton}, the map 
$\phi : \Rn \longrightarrow \inte(P)$ defined by 
\[ \phi(x)\!:=\!\!\!\!\!\!\!\!\!\left. \sum\limits_{p \text{ a vertex of } P} 
\!\!\!\!\!\!pe^{p\cdot x} \right/ \!\!\sum\limits_{q \text{ a vertex of } P} 
\!\!\!\!\!\!e^{q\cdot x} \] 
is a real analytic diffeomorphism. Composing coordinate-wise with 
the logarithm function, we then obtain that 
\[ \psi_P(x)\!:=\!\!\!\!\!\!\!\!\!\left. \sum\limits_{p \text{ a vertex of } P} 
\!\!\!\!\!\!p x^p \right/ \!\! \sum\limits_{q \text{ a vertex of } P} 
\!\!\!\!\!\! x^q \] 
yields our desired real analytic diffeomorphism from $\Rn_+$ to 
$\inte(P)$. 

The remainder of the lemma follows easily via a monomial change 
of variables. In particular, the special case where $P$ can be 
defined by a finite set of inequalities with rational coefficients 
is already embedded in the theory of toric varieties,  
e.g., \cite[Prop., Pg.\ 81]{tfulton}. The general case of 
arbitrary polytopes in $\Rn$ can be proved as follows: 
Let $d\!:=\!\dim P^w$, let $v_1,\ldots,v_{n-d}\!\in\!\Q^n$ be any linearly 
independent normal vectors of $P^w$, and let $v_{n-d+1},\ldots,v_n\!\in\!\Q^n$ 
be any linearly independent vectors parallel to $P^w$. Then, letting $A$ be 
the inverse of the $n\times n$ matrix whose $i^\thth$ column is $v_i$  
for all $i\!\in\!\{1,\ldots,n\}$, 
we can clearly write $f(y^A)\!=\!g(y_1,\ldots,y_n)$, where 
\[g(y_1,\ldots,y_n)\!:=\!\!\!\!\!\!\!\!\!
\sum_{\alpha:=(\alpha_1,\ldots,\alpha_{n-d})}\!\!\!\!\!\!\!\! 
y^{\alpha_1}_1\cdots y^{\alpha_{n-d}}_{n-d}
g_\alpha(y_{n-d+1},\ldots,y_n),\] 
the sum ranges over $\left\{(v_1\cdot a,\ldots,v_{n-d}\cdot a) \; | \; 
a\!\in\!\supp(f)\right\}$, and for any such $\alpha$ there is an 
$m_\alpha$ such that $g_\alpha$ is a $d$-variate $m_\alpha$-nomial. 
Most importantly, if 
\[\beta\!=\!\left(\min_{a\in\supp(f)} 
\{v_1\cdot a\},\ldots,\min_{a\in\supp(f)} \{v_{n-d}\cdot a\}\right)\] 
then $(y_1,\ldots,y_{n-d})^\beta 
g_\beta(y_{n-d+1},\ldots,y_n)\!=\!\init_w(f)(y^A)$. 
Clearly then, $\psi_P(Z_+(f))$ has a limit point in $\relint(P^w) 
\Longrightarrow$  
there is an $M\!>\!0$ such that $Z_+(f)$ intersects 
\[\left(\bigcap^n_{i=n-d+1} \left\{x\!\in\!\Rn_+ \; | \; 
x^{v_i}\!>\!\frac{1}{M}\right\}\right)
\cap \left(\bigcap^n_{i=n-d+1} \left\{x\!\in\!\Rn_+ \; | \; 
x^{v_i}\!<\!M\right\}\right) 
\cap \bigcap^{n-d}_{i=1}\left\{x\!\in\!\Rn_+ \; | \; 
x^{v_i}\!=\!\eps_i\right\}\] for all $\eps_i\!>\!0$ sufficiently small, 
since $P^w$ is compact and $\psi_P$ is a diffeomorphism. By Proposition 
\ref{prop:easy}, the map $x\mapsto x^A$ is a diffeomorphism, so   
there is also an $M'\!>\!0$ such that $Z_+(g)$ intersects
\[\left(\bigcap^n_{i=n-d+1} \left\{y\!\in\!\Rn_+ \; | \; 
y_i\!>\!\frac{1}{M'}\right\}\right)
\cap \left(\bigcap^n_{i=n-d+1} \left\{y\!\in\!\Rn_+ \; | 
\; y_i\!<\!M'\right\}\right)
\cap \{y_i\!=\!\delta_i \; | \; i\!\in\!\{1,\ldots,n-d\}\}\] for all 
$\delta_i\!>\!0$ sufficiently small. By our formula relating 
$g$ and $g_\beta$ (and the fact that $Z_+(g)$ is locally closed, 
being the zero set of a continuous function), we then have that  
$Z_+(g_\beta)$ intersects the last set as well. 
So $Z_+(g_\beta)$, and thus $Z_+(\init_w(f))$, is non-empty. 

To conclude, a routine monomial change of variables shows that 
\[ \psi_w(x)\!:=\left(\left. \sum\limits_{p \text{ a vertex of } 
P^w} \!\!\!\!p x^p \right/ \sum\limits_{q \text{ a vertex of } P^w} 
\!\!\!\! x^q\right)\times (x^{v_1},\ldots,x^{v_{n-d}}) \] 
gives us our desired real analytic diffeomorphism. \qed 

\noindent
Note that the converse of Lemma \ref{lemma:moment} need not hold: A simple 
counter-example is \[f(x,y)=(x^2+y^2-1)^2+(x-1)^2 \text{ and } 
w=(0,1).\]  
We also point out that the easiest way to understand the above lemma 
is to take any example $f$ with Newton polytope identical (near 
the origin) to the nonnegative orthant, and then note that 
one is in essence compactifying $\Rn_+$ by adding coordinate 
subspaces, as well as some other pieces which are images of 
$(\R^*)^k$ under monomial maps. Indeed, the monomial change of 
variables in our proof essentially results in an invertible 
affine map which sends a $d$-dimensional face of $P$ to a 
$d$-dimensional coordinate subspace of $\Rn$. 

Theorem \ref{thm:moment} then follows easily from a 
refinement of the last lemma. 
\begin{lemma} 
\label{lemma:moment2} 
Following the notation of Lemma \ref{lemma:moment}, 
assume that $Z_+(\init_w(f))$ is smooth for all 
$w\!\in\!\Rn\!\setminus\!\{\bO\}$. Then 
\begin{enumerate} 
\item{ For any facet $Q$ of $P$, every connected component of 
$\relint(Q)\cap\overline{\psi_P(Z_+(f))}$ is an $(n-2)$-manifold which is the 
set of limit points in $\relint(Q)$ of $\psi_P(C)$ for some {\bf unique} 
non-compact connected component $C$ of $Z_+(f)$. } 
\item{ $C$ a non-compact connected component of $Z_+(f) \Longrightarrow 
\psi_P(C)$ has a limit point in $\relint(Q)$ for some 
inner {\bf facet} $Q$ of $P$.} 
\end{enumerate} 
\end{lemma} 

\noindent
{\bf Proof:} To prove (1), first note that 
the last portion of Lemma \ref{lemma:moment} already tells us that 
every connected component of $\relint(Q)\cap\overline{\psi_P(Z_+(f))}$ is an 
$(n-2)$-manifold, since $Z_+(\init_w(f))$ is smooth for all $w$. (Indeed,  
the number of connected components of any $Z_+(\init_w(f))$, and thus 
$\relint(Q)\cap\overline{\psi_P(Z_+(f))}$, is finite by Lemma 
\ref{lemma:finite}.) Furthermore, it is clear 
that every connected component of $\relint(Q)\cap\overline{\psi_P(Z_+(f))}$ 
must be the set of limit points of some collection of non-compact 
components of $Z_+(f)$.  

To see why a component of $\relint(Q)\cap\overline{\psi_P(Z_+(f))}$ can be the 
limit set of just one non-compact component of $\psi_P(Z_+(f))$, we can 
specialize the monomial change of coordinates from the proof of our last 
lemma as follows: Let $w$ be any nonzero inner facet normal vector of 
$Q$ and let $A$ be any invertible $n\times n$ matrix such that $Aw$ is the 
first standard basis vector. Also let $\delta$ be the minimum value of $w\cdot 
a$ as $a$ ranges over $\supp(f)$. Lemma \ref{lemma:moment} and Proposition 
\ref{prop:easy} then tell us that $\psi_P(C)$ 
is diffeomorphic to some non-compact component of $Z_+(g)$ where 
$g(y_1,\ldots,y_n)\!=\!f(y^A)=\!y^\delta_1g_\delta(y_2,\ldots,y_n)+\sum_\alpha 
y^\alpha_1 g_\alpha(y_2,\ldots,y_n)$, the sum ranges over $\{w\cdot a\!>\!
\delta \; | \; a\!\in\!A\}$, and $y^\delta_1g_\delta(y_2,\ldots,y_n)\!=\!
\init_w(f)(y^A)$. Dividing out by $y^\delta_1$, we then obtain 
by Proposition \ref{prop:easy}, the implicit function theorem 
\cite[Thm.\ 9.28, Pg.\ 224]{rudin}, and a simple induction on 
$\dim Q$ that every connected component of $U\cap 
\overline{\psi_P(Z_+(f))}$ 
is a connected $(n-1)$-dimensional {\bf quasifold} \cite{prato,nonrat}, for 
some neighborhood $U$ of $Q$ in $P$. 
So assertion (1) is proved, with the additional strengthening that  
for all $w\!\in\!\Rn\!\setminus\!\{\bO\}$, 
every component of $\relint(P^w)\cap\overline{\psi_P(Z_+(\init_{w'}(f)))}$ is 
the limit set of some unique non-compact component of 
$Z_+(\init_{w'}(f))$, where $P^{w'}$ is a face of dimension $1+\dim P^w$ and 
$w'\!\in\!\Rn$. 

To prove (2), note that $\psi_P(C)$ must be a 
non-compact subset of $P$ and a closed subset of $\inte(P)$. 
Since $P$ is compact, $\overline{\psi_P(C)}$ must therefore be 
compact and contain a point in $\partial P$. 
So now let $Q$ be the face of highest dimension $d$ such that 
$\overline{\psi_P(C)}$ intersects $\relint(Q)$. By assertion 
(1) (and the definition of a quasifold \cite[Sec.\ 1]{prato}), $\partial 
P\cap\overline{\psi_P(C)}$ 
must be an $(n-2)$-dimensional quasifold with only finitely many connected 
components. Lemma \ref{lemma:moment} then tells us that 
$d\!<\!n-1\Longrightarrow \dim (Q\cap \overline{\psi_P(C)})\!<\!d$, 
since $\init_w(f)$ is not identically zero. So if $d\!<\!n-1$ we must then have
that $\partial P\cap \overline{\psi_P(C)}\!=\!\!\!\!\!\bigcup\limits_{Q' 
\text{ a face of } P} \!\!\!\!\relint(Q')\cap \overline{\psi_P(C)}$  
has dimension strictly less than $n-2$, thus contradicting assertion 
(1). So $d\!=\!n-1$ and assertion (2) is proved. \qed 

\noindent
Note that the smoothness hypothesis of Lemma \ref{lemma:moment2} 
in fact implies that every non-compact connected component of $Z_+(f)$ 
contains an $(n-1)$-dimensional manifold. Note also that the smoothness 
hypothesis (at least for $w$ that are inner facet normals) is necessary for 
assertion (1).  
\begin{ex} 
Consider $f(x,y)\!:=\!(x+y-1)(y-x+1)$. Then 
$Z_+(f)$ consists of a exactly 2 disjoint rays and 
$\init_{(0,1)}(f)\!=\!-(x-1)^2$ has a degenerate root at 
$1$. In particular, $(1,0)$ is a limit point of both the 
rays of this $Z_+(f)$. \dia 
\end{ex} 

\noindent
{\bf Proof of Theorem \ref{thm:moment}:} Let $Z\!:=\!Z_+(f)$. By 
Lemma \ref{lemma:moment2}, the number of non-compact connected components  
of $Z$ is no more than $\sum_w N'_w$ where the sum ranges over 
all unit inner facet normals of $P\!=\!\newt(f)$ and $N'_w$ is 
the number of connected components of $\relint(P^w)\cap 
\overline{\psi_P(Z_+(f))}$. By Lemma \ref{lemma:moment}, 
$\sum_w N'_w\!=\!\sum_w N_w$, so the first part of theorem 
\ref{thm:moment} is proved. 

The second assertion is then a trivial consequence of the first 
via the definition of $P(n,m)$. 

To conclude, assertion (0) of Theorem \ref{thm:cool} easily implies that 
for $n\!\leq\!2$ our penultimate bound specializes to exactly the number of 
points of $\supp(f)$ on the boundary of $\newt(f)$, regardless of whether $Z$ 
is smooth or not. To halve this bound, simply note that for smooth $Z$, 
every non-compact component $C$ of $Z$ is homeomorphic to an open interval. 
Therefore, by Lemma \ref{lemma:moment2}, 
$\psi_P(C)$ must intersect the boundary of $P$ exactly twice. So we are done. 
\qed 

\begin{rem} 
Bertrand Haas has pointed out that the very last assertion of
Theorem \ref{thm:moment} (concerning non-compact connected components 
of $m$-nomial curves in $\R^2_+$), in the case of {\bf integral} exponents,
follows easily from work of Isaac Newton published in 1744 
\cite[Book I, Chap.\ 3]{newton}.    
The relevant result of Newton relates Puiseux series and diagrams involving 
the portion of $\newt(f)$ visible from the origin, and can also be 
found in \cite[Chap.\ II, Paragr.\ 1, Pg.\ 213]{coolidge}. \dia 
\end{rem} 

\section*{Acknowledgements} 
The authors thank Alicia Dickenstein and Bernd Sturmfels for pointing 
out Haas' counter-example. We also thank Bertrand Haas 
for pointing out an error in an earlier version of Lemma 
\ref{lemma:imp} and supplying a nice pointer to Isaac Newton, the anonymous 
referees for giving many nice corrections, and Felipe Cucker, Jesus 
Deloera, Paulo Lima-Filho, Daniel Perrucci, and Steve Smale for some nice 
conversations. We are also indebted to Fernando Lopez Garcia, Teresa Krick, 
Daniel Perrucci, and Juan Sabia for a detailed reading of earlier versions of 
this paper. Finally, we would like thank Dima Yu.\ Grigoriev and Askold 
Georgievich Khovanski for imparting the following information about 
Konstantin Alexandrovich Sevast'yanov \cite{dima,askold}: 

\newpage 

\mbox{}
\vspace{-.3cm}
\begin{wrapfigure}{r}{1.35in}
\epsfig{file=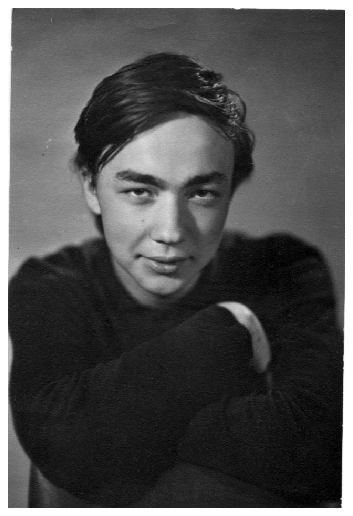,height=2in}
\end{wrapfigure} 
\mbox{}\\
\mbox{}\hspace{.5cm}Konstantin Alexandrovich Sevast'yanov was born on January 
21, 1956 in Astrakhan (an old Russian town on the Volga river) and graduated 
from a very famous mathematical high school organized by Andrey Nikolaevich 
Kolmogorov.  At the age of 17 he was a winner of the International 
Mathematical Olympiad for high school students 
% IMO? or local? 
and thus skipped his entrance exams to become a mathematics student at Moscow 
State University. His supervisor was Anatoly Georievich Kushnirenko, and 
Vladimir Igorevich Arnold and Askold Georgevich Khovanski also supervised 
Sevast'yanov's research. Sevast'yanov was a gifted student but suffered 
from poor health throughout his life. 
He formulated, around 1979, the key result that inspired Khovanski to 
create Fewnomial Theory. Sevast'yanov eventually went 
on permanent leave as his illness worsened and on December 7, 1984 he was 
killed after apparently being struck by a car. 

\vspace{-.2cm}
\noindent
\mbox{}\hspace{4.58in}\scalebox{.48}[1]{{\small Konstantin Alexandrovich 
Sevast'yanov, around 1983}}

\scalebox{.968}[1]{Those who bless us with beautiful results should never be 
forgotten, even if tragedy obscures their}\\ 
\scalebox{1.027}[1]{accomplishments. We therefore dedicate this paper 
to the memory of Konstantin Alexandrovich}\\ 
Sevast'yanov. 

\footnotesize
\bibliographystyle{acm}

\begin{thebibliography}{A}

\bibitem[Bas99]{basu} Basu, Saugata, {\it
``On Bounding the Betti Numbers and Computing the Euler
Characteristic of Semi-Algebraic Sets,''} Journal of Discrete and 
Computational Geometry, 22:1-18, (1999).

\bibitem[BP02]{nonrat} Battaglia, Fiammetta and Prato, Elisa,
{\it ``Nonrational, Nonsimple Convex Polytopes in Symplectic Geometry,''}
Math ArXiV preprint {\tt math.SG/0206149}.

\bibitem[BLR91]{blr} Benedetti, Riccardo; Loeser, Fran\c{c}ois; Risler, 
Jean-Jacques, 
{\it ``Bounding the Number of Connected Components of a Real Algebraic
Set,''} Discrete and Computational Geometry, 6:191--209 (1991).

\bibitem[BKK76]{bkk} Bernstein, David Naumovich; Kushnirenko, Anatoly 
Georievich; and Khovanski, Askold Georgevich, {\it ``Newton Polyhedra,''} 
Uspehi Mat.\ Nauk 31 (1976), no.\ 3 (189), pp.\ 201--202.

\bibitem[BCSS98]{bcss} Blum, Lenore; Cucker, Felipe; Shub, Mike; and Smale, 
Steve, {\it Complexity and Real Computation,} Springer-Verlag, 1998.     

\bibitem[BC76]{bocook} Borodin, Allan and Cook, Stephen A., {\it
``On the Number of Additions to Compute Specific Polynomials,''}
SIAM J.\ Comput.\ {\bf 5} (1976), no.\ 1, pp.\ 146--157.

\bibitem[BZ88]{buza} Burago, Yu. D. and Zalgaller, V. A., {\it
Geometric Inequalities,} Grundlehren der mathematischen Wissenschaften 285,
Springer-Verlag (1988).

\bibitem[Coo59]{coolidge} Coolidge, Julian Lowell, {\it 
A Treatise on Algebraic Plane Curves,} Dover Publications, Inc., New York 
1959. 

\bibitem[FH95]{keith} Forsythe, Keith and Hatke, Gary, {\it ``A Polynomial
Rooting Algorithm for Direction Finding,''} preprint, MIT Lincoln
Laboratories, 1995.

\bibitem[Ful84]{ifulton} Fulton, William, {\it 
Intersection Theory,} 1$^\st$ ed., Ergebnisse der Mathematik und 
ihrer Grenzgebiete 3, {\bf 2}, Springer-Verlag, 1984. 

\bibitem[Ful93]{tfulton} \underline{\hspace{\fw}}, {\it Introduction to Toric
Varieties}, Annals of Mathematics Studies, no.\ 131, Princeton University
Press, Princeton, New Jersey, 1993.   

\bibitem[GH99]{gaterhub} Gatermann, Karin and Huber, Birk, {\it 
``A Family of Sparse Polynomial Systems Arising in Chemical Reaction 
Systems,''} Preprint ZIB (Konrad-Zuse-Zentrum f\"ur Informationstechnik 
Berlin) SC-99 27, 1999. 

\bibitem[Gri82]{grigoriev} Grigor'ev, Dima Yu., {\it
``Lower Bounds in the Algebraic Complexity of Computations,''}
The Theory of the Complexity of Computations, I;
Zap.\ Nauchn.\ Sem.\ Leningrad.\ Otdel.\ Mat.\ Inst.\ Steklov (LOMI) {\bf
118} (1982), pp.\ 25--82, 214.

\bibitem[Gri00]{dima} \underline{\hspace{\dima}}, {\it
personal communication,} to J.\ Maurice Rojas, at a conference 
on Model Theoretical Algebra (Edinburgh, Scotland), Sept.\ 8, 2000. 

\bibitem[Haa02]{haas} Haas, Bertrand, {\it ``A Simple Counter-Example 
to Kushnirenko's Conjecture,''} Beitr\"age zur Algebra und Geometrie,  
Vol.\ 43, No.\ 1, pp.\ 1--8 (2002). 

\bibitem[Hir94]{hirsch} Hirsch, Morris, {\it Differential Topology,}
corrected reprint of the 1976 original, Graduate Texts in Mathematics,
33, Springer-Verlag, New York, 1994. 

\bibitem[IR96]{ir} Itenberg, Ilia and Roy, Marie-Fran\c{c}oise, {\it 
``Multivariate Descartes' Rule,''} 
Beitr\"age Algebra Geom.\ {\bf 37} (1996), no.\ 2, pp.\ 337--346. 

\bibitem[Kaz81]{kaza} Kazarnovski{\u\i}, B.\ Ja., {\it ``On Zeros of 
Exponential Sums,''} Soviet Math.\ Doklady, {\bf 23} (1981), no.\ 2, pp.\ 
347--351. 

\bibitem[Kho80]{kho} Khovanski, Askold Georgievich, {\it ``On a Class of 
Systems of Transcendental Equations,''} Dokl.\ 
Akad.\ Nauk SSSR {\bf 255} (1980), no.\ 4, pp.\ 804--807; 
English transl.\ in Soviet Math.\ Dokl.\ {\bf 22} (1980), 
no.\ 3. 

\bibitem[Kho91]{few} \underline{\hspace{\khov}}, {\it Fewnomials,}
AMS Press, Providence, Rhode Island, 1991.

\bibitem[Kho02]{askold} \underline{\hspace{\khov}}, {\it personal 
communication,} via e-mails to J.\ Maurice Rojas, August 23 -- September 1, 
2002.

\bibitem[LR97]{lr} Lagarias, Jeffrey C.\ and Richardson, Thomas J.,  
{\it ``Multivariate Descartes Rule of Signs and Sturmfels' Challenge 
Problem,''} Math.\ Intelligencer 19 (1997), no.\ 3, pp.\ 9--15. 

\bibitem[LW98]{liwang} Li, Tien-Yien and Wang, Xiaoshen, {\it ``On 
Multivariate Descartes' Rule --- A Counter-Example,''} 
Beitr\"age Algebra Geom.\ {\bf 39} (1998), no.\ 1, pp.\ 1--5.

\bibitem[Mil63]{morse} Milnor, John, {\it 
Morse Theory}, Based on lecture notes by M.\ Spivak and R.\ Wells, 
Annals of Mathematics Studies, No.\ 51 Princeton University Press, 
Princeton, N.J.\ 1963. 

\bibitem[Mil64]{milnor} \underline{\hspace{\milnor}}, 
{\it ``On the Betti Numbers of Real Varieties,''}
Proceedings of the Amer.\ Math.\ Soc.\ 15, pp.\ 275--280, 1964. 

\bibitem[Nap01]{napo} Napoletani, Domenico {\it
``A Power Function Approach to Kouchnirenko's Conjecture,''}
Contemporary Mathematics, vol.\ 286, AMS-IMS-SIAM Joint Summer Research 
Conference 
Proceedings of ``Symbolic Computation: Solving Equations in Algebra, Geometry, 
and Engineering (June 11-15, 2000, Mount Holyoke College),'' edited by 
E.\ Green, S.\ Ho\c{s}ten, R.\ Laubenbacher and V.\ Powers, AMS Press, 2001. 

\bibitem[New44]{newton} Newton, Sir Isaac, {\it Opuscula Mathematica, 
Philosophica et Philologica}, edited by G.\ F.\ Salvemini, 3 volumes, 
Lausanne and Geneva, Marc-Michel Bosquet, 1744. 

\bibitem[OP49]{op} Oleinik, O.\ and Petrovsky, I., {\it ``On the 
Topology of Real Algebraic Hypersurfaces,''} Izv.\ Akad.\ Nauk SSSR Ser.\ 
Math.\ {\bf 13} (1949), pp.\ 389--402; English transl., 
Amer.\ Math.\ Soc.\ Transl.\ (1) {\bf 7} (1962), pp.\ 399--417. 

\bibitem[Pra01]{prato} Prato, Elisa, {\it ``Simple Non-Rational Convex 
Polytopes via Symplectic Geometry,''} Topology 40 (2001), no.\ 5, 
pp.\ 961--975.  

\bibitem[Ren89]{renegar} Renegar, Jim, {\it ``On the Worst Case
Arithmetic Complexity of Approximating Zeros of Systems of Polynomials,''}
SIAM J.\ Comput.\ {\bf 18} (1989), no.\ 2, pp.\ 350--370.

\bibitem[Ris85]{risler} Risler, Jean-Jacques, {\it ``Additive Complexity
and Zeros of Real Polynomials,''} SIAM J.\ Comput.\ {\bf 14} (1985), no.\ 1,
pp.\ 178--183.

\bibitem[Roj99]{jpaa} Rojas, J.\ Maurice, {\it ``Toric
Intersection Theory for Affine Root Counting,''} Journal of Pure and
Applied Algebra, vol.\ 136, no.\ 1, March, 1999, pp.\ 67--100.

\bibitem[Roj00a]{real} \underline{\hspace{\jmr}}, {\it ``Some Speed-Ups and 
Speed Limits for Real Algebraic Geometry,''} Journal of Complexity, FoCM 1999 
special issue, vol.\ 16, no.\ 3 (sept.\ 2000), pp.\ 552--571. 

\bibitem[Roj01]{myadic} \underline{\hspace{\jmr}}, {\it ``Finiteness 
for Arithmetic Fewnomial Systems,''} invited paper, Contemporary Mathematics, 
vol.\ 286, 
AMS-IMS-SIAM Joint Summer Research Conference Proceedings of ``Symbolic 
Computation: Solving Equations in Algebra, Geometry, and Engineering (June 
11--15, 2000, Mount Holyoke College),'' edited by E.\ Green, S.\ Ho\c{s}ten, 
R.\ Laubenbacher and V.\ Powers, AMS Press, 2001. 

\bibitem[Roj02]{ari} \underline{\hspace{\jmr}}, {\it ``Arithmetic 
Multivariate Descartes' Rule,''} Math ArXiV preprint {\tt math.NT/0110327}, 
submitted for publication. 

\bibitem[RY02]{rojasye} Rojas, J.\ Maurice and Ye, Yinyu, {\it 
``On Solving Fewnomials Over an Interval in Fewnomial Time,''} 
Math ArXiV preprint {\tt math.NA/0106225}, submitted for publication. 

\bibitem[Rud76]{rudin} Rudin, Walter, {\it Principles of   
Mathematical Analysis,} 3$^\rd$ edition, McGraw-Hill, 1976.  

\bibitem[Sha77]{sha} Shafarevich, Igor R., {\it Basic Algebraic 
Geometry,} Springer Study Edition, Springer-Verlag, 1977. 

\bibitem[Sma70]{smale} Smale, Steve, {\it ``Topology and Mechanics I,''} 
Invent.\ Math.\ {\bf 10} (1970), pp.\ 305--331. 

\bibitem[SL54]{descartes} Smith, David Eugene and Latham, Marcia L., {\it 
The Geometry of Ren\'e Descartes,} translated from the French and Latin 
(with a facsimile of Descartes' 1637 French edition), 
Dover Publications Inc., New York (1954).  

\bibitem[Sou70]{souriau} Souriau, J.-M., {\it Structure des Syst\`emes 
Dynamiques,} Dunod, Paris, 1970, Ma$\hat{\text{\i}}$trises de Math\'ematiques. 

\bibitem[Stu94]{berndviro} Sturmfels, Bernd, {\it ``On the Number of Real Roots 
of a Sparse Polynomial System,''} Hamiltonian and Gradient Flows, Algorithms 
and Control, pp.\ 137--143, Fields Inst.\ Commun., 3, Amer.\ Math.\ Soc., 
Providence, RI, 1994. 
 
\bibitem[Stu98]{poly} \underline{\hspace{\bernd}}, {\it ``Polynomial Equations 
and Convex Polytopes,''} American Mathematical Monthly {\bf 105} (1998), 
no.\ 10, pp.\ 907--922.

\bibitem[Tho65]{thom} Thom, Ren\'e, {\it ``Sur l'homologie 
des vari\'et\'es alg\'ebriques r\'eelles,''} In S.\ 
Cairns (Ed.), Differential and Combinatorial Topology, 
Princeton University Press, 1965. 

\bibitem[Ver99]{verschelde} Verschelde, Jan, {\it ``Algorithm 795: PHCpack: A 
General-Purpose Solver for Polynomial Systems by Homotopy Continuation,''} 
ACM Transactions on Mathematical Software 25(2): pp.\ 251-276, 1999. 

\bibitem[VR02]{sagarrojas} Vidyasagar, M.\ and Rojas, J.\ Maurice, 
{\it ``An Improved Bound on the VC-Dimension of Neural Networks with 
Polynomial Activation Functions,''} Math ArXiV preprint 
{\tt math.OC/0112208}, submitted for publication.   

\bibitem[Vir84]{viro} Viro, Oleg Ya., {\it ``Gluing of Plane Real Algebraic 
Curves and Constructions of Curves of Degrees $6$ and $7$,''} 
Topology (Leningrad, 1982), pp.\ 187--200, Lecture Notes in Math., 1060, 
Springer, Berlin, 1984. 

\bibitem[Zel99]{zell} Zell, Thierry, 
{\it ``Betti Numbers of Semi-Pfaffian Sets,''} 
Effective Methods in Algebraic Geometry (MEGA '98, Saint-Malo, 1998). 
J.\ Pure Appl.\ Algebra {\bf 139} (1999), no.\ 1--3, pp.\ 323--338. 

\end{thebibliography}

\end{document}